\documentclass{amsart}
\usepackage{amssymb,amsfonts}
\usepackage{enumerate}
\usepackage{mathrsfs}
\usepackage{graphicx}
\usepackage{mathtools}
\usepackage{leftidx}
\usepackage{tikz-cd}
\usepackage{bm}
\usepackage{url}
\newtheorem{thm}{Theorem}[section]
\newtheorem{cor}[thm]{Corollary}
\newtheorem{prop}[thm]{Proposition}
\newtheorem{lem}[thm]{Lemma}

\theoremstyle{definition}
\newtheorem{defn}[thm]{Definition}

\newtheorem{exmp}[thm]{Example}

\theoremstyle{remark}
\newtheorem{rem}[thm]{Remark}

\makeatletter
\let\c@equation\c@thm
\makeatother
\numberwithin{equation}{section}

\def\bthm{\begin{thm}}
\def\ethm{\end{thm}}
\def\blm{\begin{lem}}
\def\elm{\end{lem}}
\def\bdf{\begin{defn}}
\def\edf{\end{defn}}
\def\bpf{\begin{proof}}
\def\epf{\end{proof}}
\def\bpp{\begin{prop}}
\def\epp{\end{prop}}
\def\bcor{\begin{cor}}
\def\ecor{\end{cor}}
\def\brm{\begin{rem}}
\def\erm{\end{rem}}
\def\beg{\begin{exmp}}
\def\eeg{\end{exmp}}
\def\bD{\mathbb{D}}

\def\bL{\mathbb{L}}
\def\bM{\mathbb{M}}
\def\bN{\mathbb{N}}

\def\bQ{\mathbb{Q}}

\def\bZ{\mathbb{Z}}
\def\cA{\mathcal{A}}
\def\cB{\mathcal{B}}

\def\cD{\mathcal{D}}

\def\cF{\mathcal{F}}

\def\cM{\mathcal{M}}

\def\cO{\mathcal{O}}
\def\cP{\mathcal{P}}
\def\cQ{\mathcal{Q}}
\def\cR{\mathcal{R}}
\def\cS{\mathcal{S}}

\def\scC{\mathscr{C}}

\def\scO{\mathscr{O}}
\def\scP{\mathscr{P}}

\def\scX{\mathscr{X}}

\def\frf{\mathfrak{f}}
\def\frg{\mathfrak{g}}

\newcommand{\raq}{\,\rightarrow \,}

\newcommand{\rintoq}{\,\hookrightarrow\,}

\newcommand{\xraq}[2][]{\, \xrightarrow[#1]{#2} \,}

\newcommand{\xlintoq}[2][]{\, \xhookleftarrow[#1]{#2} \,}

\newcommand{\ra}{\rightarrow}
\newcommand{\la}{\leftarrow}
\newcommand{\rinto}{\hookrightarrow}

\newcommand{\ronto}{\twoheadrightarrow}

\newcommand{\rsa}{\stackrel{\sim}{\rightarrow}}

\newcommand{\xra}[2][]{\xrightarrow[#1]{#2}}

\newcommand{\xronto}[2][]{\xrightarrow[#1]{#2}\mathrel{\mkern-14mu}\rightarrow}

\newcommand{\Set}{{\rm Set}}
\newcommand{\Cat}{{\rm Cat}}
\newcommand{\Setdel}{{\rm Set}_{\Delta}}
\newcommand{\Mod}{{\rm Mod}}
\newcommand{\Vect}{{\rm Vect}}
\newcommand{\Alg}{{\rm Alg}}
\newcommand{\CAlg}{{\rm CAlg}}

\newcommand{\Sym}{{\rm Sym}}

\newcommand{\Tr}{{\rm Tr}}

\newcommand{\Ch}{{\rm Ch}}
\newcommand{\DGA}{{\rm dga}}

\newcommand{\CDGAn}{{\rm cdga}^{\leq 0}}
\newcommand{\dgCat}{{\rm dgcat}}

\newcommand{\dgCatn}{{\rm dgcat}^{\leq 0}}

\newcommand{\dgCatnOk}{{\rm dgcat}^{\leq 0}_{\scO,k}}

\newcommand{\Moddg}{\underline{{\rm Mod}}}

\newcommand{\dgcat}{{\rm dgcat}}

\newcommand{\cone}{{\rm cone}}

\newcommand{\RHomcom}{{\bm R} \underline{{\rm Hom}}}

\newcommand{\CN}{CN}

\newcommand{\Ob}{{\rm Ob}}
\newcommand{\op}{{\rm op}}

\newcommand{\id}{{\rm id}}

\newcommand{\Hom}{{\rm Hom}}
\newcommand{\Homcom}{\underline{{\rm Hom}}}
\newcommand{\End}{{\rm End}}
\newcommand{\Endcom}{\underline{{\rm End}}}

\newcommand{\Fun}{{\rm Fun}}

\newcommand{\colim}{{\rm colim}}

\newcommand{\cHom}{\mathscr{H}\text{\kern -3pt {\calligra\large om}}\,}

\newcommand{\Map}{{\rm Map}}

\newcommand{\Ho}{{\rm Ho}}

\newcommand{\Spec}{{\rm Spec}}

\newcommand{\GL}{{\rm GL}}

\newcommand{\Rep}{{\rm Rep}}

\newcommand{\DRep}{{\rm DRep}}

\newcommand{\dAff}{{\rm dAff}}

\newcommand{\QCoh}{{\rm QCoh}}

\newcommand{\Chdg}{\underline{\Ch}}

\newcommand{\ab}{{\rm ab}}

\newcommand{\cAe}{\cA^{e}}

\newcommand{\cSA}{\cS(\cA)}

\newcommand{\MD}{\bD}

\newcommand{\cyc}{{\rm cyc}}

\newcommand{\tot}{{\rm tot}}
\newcommand{\dtot}{d_{{\rm tot}}}

\newcommand{\grMod}{{\rm grMod}}

\newcommand{\cSAfr}{\cS_{\scO'\text{-fr}}(\cA)}
\newcommand{\Opfr}{\scO'\text{-fr}}
\newcommand{\fr}{\text{fr}}

\newcommand{\com}{\text{com}}

\newcommand{\Car}{\text{Car}}

\newcommand{\DR}{\text{DR}}

\newcommand{\cl}{\text{cl}}

\newcommand{\Pol}{\text{Pol}}

\DeclareMathAlphabet{\mathpzc}{OT1}{pzc}{m}{it}

\title{Pre-Calabi-Yau structures and moduli of representations}
\author{Wai-Kit Yeung}
\address{Kavli IPMU, The University of Tokyo}
\email{wai-kit.yeung@ipmu.jp}
\begin{document}

\begin{abstract}
We establish a system of formal noncommutative calculus for differential forms and polyvector fields, which forms the foundations for the study of pre-Calabi-Yau categories. Using an explicit trace map, we show that any $n$-Calabi-Yau structure on a non-positively graded dg algebra $A$ induces a $(2-n)$-shifted symplectic structure on its derived moduli stack of representations; while any $n$-pre-Calabi-Yau structure on $A$ induces a $(2-n)$-shifted Poisson structure on this derived moduli stack.
\end{abstract}

\maketitle


\tableofcontents

\section{Introduction}

We achieve three objectives in this paper:	
\begin{enumerate}
	\item[(a)] Establish a system of formal noncommutative calculus for differential forms and polyvector fields. In particular, this sets up the foundations for the study of pre-Calabi-Yau categories, a notion first appeared in \cite{KV13, KTV}, and independently studied by us.
	\item[(b)] For any associative algebra (or more generally non-positively graded dg algebra) $A$, give an explicit global quotient that represents its derived moduli stack $\mathpzc{DRep}(\cA)$ of finite dimensional representations.
	\item[(c)] Using an explicit trace map, show that any $n$-Calabi-Yau structure on $A$ induces a $(2-n)$-shifted symplectic structure on $\mathpzc{DRep}(A)$; while any $n$-pre-Calabi-Yau structure on $A$ induces a $(2-n)$-shifted Poisson structure on $\mathpzc{DRep}(A)$.
\end{enumerate}

These three objectives are interconnected with each other. Together, they form a coherent picture of formal noncommutative algebraic geometry. We briefly explain the underlying philosophy in Section \ref{sec_FNAG}.
Now we discuss these three objectives separately. 
In this introduction, we work over a fixed field $k$ of characteristic zero.

In \cite{Yeu3}, we established an extension of the Feigin-Tsygan Theorem to higher Hodge degrees. This gives a description of the negative and periodic cyclic homology of a dg algebra $A$ in terms of the extended noncommutative de Rham complex \cite{GS12} of any cofibrant resolution $Q \xra{\sim} A$.
This can be viewed as a basic result about formal noncommutative differential forms. Our objective (a) can be viewed as a further development of this picture, where we consider also the dual side of polyvector fields. 

The formal noncommutative analogues of both differential forms and polyvector fields share some features with their commutative counterparts. 
More precisely, on the commutative side, we consider a derived stack $X$, then the de Rham complex is an $\bN$-graded mixed complex
\begin{equation}  \label{DR_X_mixed_intro}
	\DR^{(0)}(X) \xraq{D} \DR^{(1)}(X) \xraq{D} \DR^{(2)}(X) \xraq{D} \ldots
\end{equation}
while the $m$-shifted polyvector fields $\Pol^{\bullet}(X,m)$ admits a bracket of weight $-1$:
\begin{equation}   \label{Pol_X_bracket_intro}
	\{ -,-\} \, : \, \Pol^{p}(X,m)[m+1] \otimes \Pol^{q}(X,m)[m+1] \raq \Pol^{p+q-1}(X,m)[m+1]
\end{equation}
making $\bigoplus_{r \geq 1} \Pol^{r}(X,m)[m+1]$ a weight graded dg Lie algebra.

For the noncommutative analogue, consider a small dg category $\cA$. Then in place of \eqref{DR_X_mixed_intro} there is an $\bN$-graded mixed complex
\begin{equation}  \label{scX_A_mixed_intro}
	\scX^{(0)}(\cA) \xraq{B} \scX^{(1)}(\cA) \xraq{B} \scX^{(2)}(\cA) \xraq{B} \ldots
\end{equation}
which is, up to shifts, the same as the extended noncommutative de Rham complex. 

Also, in place of \eqref{Pol_X_bracket_intro}, there is a weight graded dg Lie algebra $\bigoplus_{r \geq 1} \scP^{r}(\cA,m)[m+1]$, called the extended necklace dg Lie algebra, with a bracket
\begin{equation}   \label{scP_A_bracket_intro}
	\{ -,-\} \, : \, \scP^{p}(\cA,m)[m+1] \otimes \scP^{q}(\cA,m)[m+1] \raq \scP^{p+q-1}(\cA,m)[m+1]
\end{equation}

In \cite{Yeu3}, we showed that \eqref{scX_A_mixed_intro} is closely related to the negative and periodic cyclic homology of $\cA$. This allows one to characterize Calabi-Yau structures in terms of \eqref{scX_A_mixed_intro} (see Definition \ref{CY_str_def}). In a similar vein, $\scP^{p}(\cA,m)$ is closely related to the higher Hochschild cohomology of $\cA$ (see Theorem \ref{ext_necklace_Lie_subalg}). This allows one to characterize pre-Calabi-Yau structures in terms of \eqref{scP_A_bracket_intro}.

Now we discuss objective (b). For any associative algebra $A$, its moduli stack of $r$-dimensional representations can be realized as a global quotient
$\mathpzc{Rep}_r(A)= [ \Spec(A_r) / \GL(r) ]$, where $A_r$ is the commutative algebra obtained by putting the subscript $(-)_{ij}$ (for $1\leq i,j\leq r$) to every generator and relation of $A$ (see Section \ref{sec_DRep} for details). 
The procedure $A \mapsto A_r$ of putting subscripts $(-)_{ij}$ makes sense for any dg algebra $Q$, in which case $Q_r$ is a commutative dg algebra. In particular, applying this to any cofibrant resolution $Q \xra{\sim} A$, we can consider the derived stack $\mathpzc{DRep}_r(A)= [ \Spec(Q_r) / \GL(r) ]$. 
In Section \ref{DRep_subsec}, we show that this recovers the derived moduli stack of pseudo-perfect dg modules on $A$ with Tor-amplitude contained in $[0,0]$, as defined in \cite{TV07} (see Theorem \ref{DRep_moduli_thm}):
\begin{equation*}
	\cM_{A}^{[0,0]} \, \simeq \,  \mathpzc{DRep}(A) \, := \,
	\coprod_{r \geq 0} 
	\mathpzc{DRep}_r(A)
\end{equation*}

Three classes of examples of moduli spaces are widely studied: (i) moduli space of (complexes of) sheaves on a scheme; (ii) moduli space of representations of an associative algebra; (iii) moduli space of local systems on a topological space. 
Our approach to moduli spaces is directly applicable to (ii). It also covers (iii) because one can take $A$ to be the dg algebra of chains on the based loop space. More generally, by taking $A$ to be a suitable twist of the exit path category \cite{Lur}, our approach is also in some cases applicable to the moduli space of perverse sheaves with a given stratification (see \cite{BEY} for an interesting example). 
Our approach is not directly applicable to (i). However, often in local calculations, (i) is reduced to (ii) by certain quiver presentation, in which case our approach is also applicable.
While the general framework in \cite{TV07, TV08} covers all three cases and more, our approach has the advantage of being explicit.

Objective (c) concerns the main result of this paper:
\bthm[=Theorems \ref{CY_induce_thm}, \ref{PCY_induce_thm}]  \label{CY_induce_thm_intro}
Any $n$-Calabi-Yau structure on $A$ induces a $(2-n)$-shifted symplectic structure on $\mathpzc{DRep}(A)$; while any $n$-pre-Calabi-Yau structure on $A$ induces a $(2-n)$-shifted Poisson structure on $\mathpzc{DRep}(A)$.
\ethm

We now briefly discuss the proof of Theorem \ref{CY_induce_thm_intro}. Since $\mathpzc{DRep}_r(A)$ is presented as a global quotient, we first describe shifted symplectic and Poisson structures on global quotients $X = [Y/G]$, where $Y = \Spec \, B$ for a non-positively graded commutative dg algebra $B$ (we assume that $B$ is cofibrant) and $G$ is a reductive group. It is proved in \cite{Yeu4} that, in such cases, the mixed complex \eqref{DR_X_mixed_intro} and the dg Lie algebra \eqref{Pol_X_bracket_intro} can be represented by explicit ``Cartan models'', denoted by $\DR_{\Car}^{\bullet}(Y/G)$ and $\Pol_{\Car}^{\bullet}(Y/G,m)$ respectively, which we use to characterize shifted symplectic and Poisson structures. 

In the case $X = \mathpzc{DRep}_r(A) = [Y/G]$ for $Y = \DRep_r(A)$ and $G = \GL_r$, we construct explicit trace maps
\begin{equation*}
	\begin{split}
	\scX^{\bullet}(A) &\raq \DR_{\Car}^{\bullet}(Y/G)  \\
	\scP^{\bullet}(\cA,m) &\raq \Pol_{\Car}^{\bullet}(Y/G,m)
	\end{split}
\end{equation*}
which is respectively a map of $\bN$-graded mixed complex and a map of weight graded dg Lie algebras.
Theorem \ref{CY_induce_thm_intro} is proved using these explicit trace maps.

Results similar to Theorem \ref{CY_induce_thm_intro} have appeared in the literature. In particular, the first statement (the one regarding shifted symplectic structures) was proved in \cite{BD21} in a more general form: it concerns the whole derived stack $\cM_A$ without restriction on Tor-amplitudes. Related results have also already appeared in \cite{PTVV}. While our version of this first statement is less general than the one in \cite{BD21}, our approach is more explicit, and is rather different from the categorical approach in \cite{BD21}. Thus, we think that our result has independent value. The second statement (the one regarding shifted Poisson structures) seems to be new. We expect it to be useful in questions of quantization.

\vspace{0.2cm}

\textbf{Acknowledgement.} The author thanks Yuri Berest, Christopher Brav, Sheel Ganatra, Ezra Getzler, Nick Rozenblyum and Boris Tsygan for helpful discussions.

\section{Some noncommutative calculus}
Throughout this paper, we fix a commutative ring $k$ with unit. Unadorned tensor product will be understood to be over $k$. Starting from Section \ref{sec_DRep}, we will assume that $k$ is a field of characteristic zero. We will work with cohomological gradings, {\it i.e.}, the differential increases the degree by $1$.
We will follow the conventions and notations of dg categories in \cite{Yeu1, Yeu3}, except that {\it loc.cit.} uses homological grading. We recall some of these conventions and notations now. Details can be found in \cite{Yeu1}.

Given a small dg category $\cA$, a module will mean a right dg module, {\it i.e.,} a dg functor $\cA^{\op} \ra \Chdg(k)$. The category of modules will be denoted as $\Mod(\cA)$. Given $M,N \in \Mod(\cA)$, then a \emph{pre-map} of (cohomological) degree $i$ is a map $f : M \ra N[i]$ of graded modules over the graded category $\cA$ ({\it i.e.,} we do not require $f$ to commute with the differentials). Then the pre-maps of various degrees form a cochain complex $\Homcom_{\cA}(M,N)$ by $d(f) = d\circ f - (-1)^{|f|} f \circ d$. Denote by $\Moddg(\cA)$ the dg category of modules, with Hom complexes given by $\Homcom_{\cA}(M,N)$. Thus, $\Mod(\cA) = Z^0(\Moddg(\cA))$.

Given a bimodule $M \in \Mod(\cAe)$, we write ${}_yM_x := M(x,y)$. Its \emph{naturalization} is the cochain complex $M_{\natural} \in \Ch(k)$ defined by
\begin{equation*}
	M_{\natural} \, := \, M \otimes_{\cAe} \cA \, = \, \Bigl( \, \bigoplus_{x \in \Ob(\cA)} M(x,x) \Bigr) / (\xi f - (-1)^{|f||\xi|} f \xi)_{f \in \cA(x,y), \xi \in M(y,x)}
\end{equation*}
where we quotient out by the $k$-linear span of the displayed relations.

Write $\scO := \Ob(\cA)$ and ${}_{\bullet}M \otimes_{\scO} N_{\bullet} := \bigoplus_{x \in \scO} {}_{\bullet}M_x \otimes {}_x N_{\bullet}$.
For $n \geq 0$, define $\cR_n(\cA) \in \Mod(\cAe)$ by
\begin{equation}  \label{RnA_def}
	\cR_n(\cA) \, := \, \cA \otimes_{\scO} \stackrel{(n+2)}{\ldots} \otimes_{\scO} \cA
\end{equation}
Then $\cR_{\bullet}(\cA)$ forms a simplicial object in $\Mod(\cAe)$. The total object of the associated complex is denoted as $\cR(\cA) \in \Mod(\cAe)$. In other words, we have
\begin{equation*}
	\cR(\cA) \, = \, \bigoplus_{n \geq 0} \, \cR_n(\cA)[n]
\end{equation*}
with a suitable differential. There is a canonical quasi-isomorphism $\cR(\cA) \xra{\sim} \cA$ in $\Mod(\cAe)$. The bimodule $\cR(\cA)$ is called the \emph{bar resolution}.

We say that a cochain complex $V$ is h-flat if $V \otimes - : \Ch(k) \ra \Ch(k)$ preserves acyclic complexes. We say that $\cA$ is $k$-flat if each $\cA(x,y) \in \Ch(k)$ is h-flat. We say that $\cA$ is \emph{linearly cofibrant} if each $\cA(x,y) \in \Ch(k)$ is cofibrant.
If $\cA$ is $k$-flat, then $\cR(\cA)$ is a flat bimodule resolution of $\cA$. If $\cA$ is linearly cofibrant, then $\cR(\cA)$ is a cofibrant bimodule resolution of $\cA$.

Define $\cSA \in \Mod(\cAe)$ to be the cone
\begin{equation}  \label{cSA_def}
	\cSA \, := \, \cone \, [ \, \Omega^1(\cA) \xraq{\alpha} \cA \otimes_{\scO} \cA  \, ] 
\end{equation}
where $\Omega^1(\cA)$ is the bimodule of noncommutative Kahler differentials, and the map $\alpha$ is defined by $\alpha(Df) = f \otimes 1_x - 1_y \otimes f$ for $f \in \cA(x,y)$. 
There is a canonical quasi-isomorphism $\cSA \xra{\sim} \cA$ in $\Mod(\cAe)$. We call $\cSA$ the \emph{short resolution}, or \emph{Cuntz-Quillen resolution} of $\cA$.
For $x \in \scO$, we will write $E_x := 1_x \otimes 1_x \in {}_x(\cA \otimes_{\scO} \cA)_x$. 

We say that $\cA$ is \emph{almost cofibrant} if it is linearly cofibrant and if $\cSA \in \Mod(\cAe)$ is cofibrant. Every cofibrant dg category is almost cofibrant (see, {\it e.g.}, \cite{Yeu1}).

If $\cA$ is $k$-flat, then it is said to be \emph{smooth} if $\cA$ is perfect as a bimodule over itself. In general, $\cA$ is said to be \emph{smooth} if some, and hence any, $k$-flat $\cA' \in \dgcat_k$ quasi-equivalent to it is smooth.

\subsection{Formal noncommutative algebraic geometry}  \label{sec_FNAG}

In this subsection, we discuss some general principles of formal noncommutative algebraic geometry. These are useful organizational principles that has informed the development in this paper and elsewhere, but since these are of a semi-philosophical nature, we will keep the discussion as brief as possible.

Formal noncommutative algebraic geometry has two aspects. The first aspect may be called the \emph{ontological aspect}, which, roughly speaking, says that one should formally replace commutative algebra/shemes/derived stacks by dg categories/$A_{\infty}$-categories, and quasi-coherent sheaves by bimodules, and develop formal noncommutative analogues of usual commutative notions. The second aspect may be called the \emph{phenomenological aspect}, encapsulated by the Kontsevich-Rosenberg principle, which says that a noncommutative analogue of a property/structure $\cP$ should be a property/structure $\cP_{{\rm nc}}$ on associative algebras (or more generally dg categories/$A_{\infty}$-categories) $A$ that would induce the property/structure $\cP$ on the moduli space of representations of $A$.
The strength of the Kontsevich-Rosenberg principle is that these two aspects end up doing the same thing.

We will be developing noncommutative differential calculus in this section. Recall the short resolution $\cSA$ from \eqref{cSA_def}. The starting point will be the following principle:
\begin{equation} \label{NC_prin_cSA}
	\parbox{40em}{The bimodule $\cSA[-1]$ is a noncommutative analogue of the sheaf of Kahler differentials $\Omega^1_{{\rm com}}(X)$ on a scheme $X$.}
\end{equation}

This can be justified by both aspects of formal noncommutative algebraic geometry as discussed above. For the phenomenological aspect, see Section \ref{sec_sympl_poiss} below. For the ontological aspect, let us simply mention that \eqref{NC_prin_cSA} is an aesthetically pleasing choice in view of examples of: (i) dg algebras arising from Koszul duality, more precisely the cobar construction; and (ii) chain dg algebras of based loop spaces. 

Another useful principle is 
\begin{equation}  \label{NC_prin_naturalize}
	\parbox{40em}{The naturalization procedure $M \mapsto M_{\natural}$ of a bimodule $M$ is the noncommutative analogue of the procedure $\cF \mapsto \Gamma(X;\cF)$ of taking global sections.}
\end{equation}

Again, this  can be justified by both aspects of formal noncommutative algebraic geometry. For the phenomenological aspect, see Section \ref{subsec_trace} below. For the ontological aspect, recall that if $N$ is a finitely generated projective module over a commutative ring $B$, then the total space of the corresponding vector bundle is the $\Spec$ of $\Sym_B(N^{\vee})$, so that the set of sections is $\Hom_{B \downarrow \CAlg_k}( \Sym_B(N^{\vee}) , B ) = (N^{\vee})^{\vee} = N$.
In the noncommutative case, if $A$ is an associative algebra, and $M$ is an ordinary ({\it i.e.,} non-dg) bimodule that is finitely generated and projective, then the total section may be replaced by
$\Hom_{A \downarrow \Alg_k}( T_{A}(M^{\vee}) , \cA ) = \Hom_{A^e}(M^{\vee} , \cA) = M_{\natural}$. 

The principles \eqref{NC_prin_cSA} and \eqref{NC_prin_naturalize} should be thought of as useful rules of thumb, instead of part of a system of axioms. 
Eventually, the main justification for \eqref{NC_prin_cSA} and \eqref{NC_prin_naturalize} is that these principles guide one towards interesting mathematical structures and results. This is indeed the main theme of the present paper.

\subsection{The extended noncommutative de Rham complex}  \label{sec_nc_DR}
In this subsection, we recall some results of \cite{Yeu3} concerning the extended noncommutative de Rham complex as well as its relations with negative and periodic cyclic homology. 

For a smooth variety, global differential forms are the same as global functions on the $1$-shifted tangent bundle. Thus, the principles \eqref{NC_prin_cSA} and \eqref{NC_prin_naturalize} lead us to consider the naturalization of $T_{\cA}(\cSA)$. There are two ways to take such a naturalization: either as a bimodule over $\cA$, for which we will write $T_{\cA}(\cSA)_{\natural}$; or as a bimodule over itself, for which we will write $T_{\cA}(\cSA)_{\cyc}$. We will mostly consider the latter (see \cite{Yeu3} for a comparison of the two), which decomposes as
\begin{equation*}
	T_{\cA}(\cSA)_{\cyc} \, = \, \bigoplus_{n \geq 0} \scX^{(n)}(\cA)
\end{equation*}
where the weight $n$ component is given by
\begin{equation}  \label{scX_n_def}
	 \scX^{(n)}(\cA) \, = \, (\cSA \otimes_{\cA} \stackrel{(n)}{\ldots} \otimes_{\cA} \cSA)_{\natural, C_n}
\end{equation}
where we first take the naturalization $(-)_{\natural}$ of the $\cA$-bimodule $\cSA \otimes_{\cA} \stackrel{(n)}{\ldots} \otimes_{\cA} \cSA$, which then inherits an action by the cyclic group $C_n := \bZ/n \bZ$, for which we take the coinvariants.

The differential in the cochain complex \eqref{scX_n_def} will be denoted as $b$. The complex $(\scX^{(n)}(\cA),b)$ may be regarded as (up to shift) a noncommutative analogue of global $n$-forms. There is an analogue of the de Rham differential of the form
\begin{equation}  \label{scX_mixed}
	\scX^{(0)}(\cA) \xraq{B} \scX^{(1)}(\cA) \xraq{B} \scX^{(2)}(\cA) \xraq{B} \ldots
\end{equation}
making the collection $(\scX^{\bullet}(\cA),b,B)$ an $\bN$-graded mixed complex in the sense of the following

\bdf
An \emph{$\bN$-graded mixed complex} is a sequence $\{ (C^{(n)},b)\}_{n \geq 0}$ of cochain complexes, together with maps $B: C^{(n)} \ra C^{(n+1)}$ of cohomological degree $-1$, satisfying $B^2 = 0$ and $Bb + bB = 0$.
\edf

The definition and basic properties of the maps $B$ in \eqref{scX_mixed} may be found in \cite{Yeu3}. Essentially, it is defined by writing the letter $D$ in front of the expression, and simplify by using the Leibniz rule. 

\bdf
Given an $\bN$-graded mixed complex $(C^{\bullet},b,B)$, its \emph{(direct product) total complex} is the complex
\begin{equation*}
	C^{\tot} \, := \, \prod_{n \geq 0} \, C^{(n)} \cdot u^n \, , \qquad \qquad \dtot = b + uB
\end{equation*}
where $u$ is a variable of (cohomological) degree $2$. It comes with a filtration
\begin{equation*}
	F^r C^{\tot} \, := \, \prod_{n \geq r} \, C^{(n)} \cdot u^n \, , \qquad \qquad \dtot = b + uB
\end{equation*}
\edf

Denote by $\CN(\cA)$ the negative cyclic complex of $\cA$.
We recall a result of \cite{Yeu3}:

\bthm  \label{X_complex_thm_2}
Suppose that $\cA$ is almost cofibrant, and that $\bQ \subset k$.
Then for each $r \geq 1$, there is a zig-zag of quasi-isomorphisms relating 
$F^r \scX^{\tot}(\cA)$ and $\CN(\cA) \cdot u^r = \CN(\cA)[-2r]$.
\ethm

\brm
The entire complex $\scX^{\tot}(\cA)$ is ``almost'' quasi-isomorphic to the periodic cyclic complex, in the sense that the reduced version of $\scX^{\tot}(\cA)$ is quasi-isomorphic to the reduced periodic cyclic complex. See \cite{Yeu3} for details.
\erm

We now formulate a noncommutative analogue of the notion of shifted symplectic structures along the lines of the principles \eqref{NC_prin_cSA} and \eqref{NC_prin_naturalize}. By Theorem \ref{X_complex_thm_2}, it will correspond to a (left) Calabi-Yau structure in the sense of \cite{BD19}. 

Write $X^{(2)}(\cA) := (\cSA \otimes_{\cA}  \cSA)_{\natural}$, so that $\scX^{(2)}(\cA) = X^{(2)}(\cA)_{C_2}$. Consider the map
\begin{equation*}
\rho = \id + \tau \, : \, \scX^{(2)}(\cA) \, = \, X^{(2)}(\cA)_{C_2} \raq X^{(2)}(\cA)^{C_2}
\end{equation*}

From this, we see that any cocycle $\omega \in Z^{-n}( \scX^{(2)}(\cA) )$ determines a map of bimodules (see \eqref{bimod_dual} below for the definition of bimodule dual)
\begin{equation}  \label{CY_omega_map}
	\rho(\omega)^{\#} \, : \, \cSA^{\vee}[n] \raq \cSA
\end{equation}

\bdf  \label{CY_str_def}
Suppose that $\cA$ is almost cofibrant and smooth, and that $\bQ \subset k$.  
Then an \emph{$n$-Calabi-Yau structure} on $\cA$ is a cocycle $\widetilde{\omega} \in Z^{4-n} (F^2 \scX^{\tot}(\cA))$ whose lowest order term $\omega \in Z^{-n}(\scX^{(2)}(\cA))$ is \emph{non-degenerate} in the sense that the induced map \eqref{CY_omega_map} is a quasi-isomorphism.
Two such $n$-Calabi-Yau structures are said to be equivalent if they represent the cohomology class $H^{4-n} (F^2 \scX^{\tot}(\cA))$.

More generally, if $\cA$ is smooth and $\bQ \subset k$, then an \emph{$n$-Calabi-Yau structure} on $\cA$ is an $n$-Calabi-Yau structure on any (almost) cofibrant resolution $\cQ \xra{\sim} \cA$. The set of equivalence classes of $n$-Calabi-Yau structures is independent of the choice of resolution.
\edf

\subsection{Multiduals}

Recall that (see \cite{Yeu1} for notations), given a bimodule $M \in \Mod(\cAe)$ over a small dg category $\cA$, its \emph{bimodule dual} is the bimodule $M^{\vee} \in \Mod(\cAe)$ given by
\begin{equation}  \label{bimod_dual}
	{}_{\bullet}(M^{\vee})_{\bullet} \, := \, \Homcom_{\cAe}( \, {}_L M_{R} \, , \, {}_{\bullet}\cA_R \otimes {}_L \cA_{\bullet} \, )
\end{equation}
where we take the Hom complex of pre-maps from $M$ to $\cA \otimes \cA$ that are $\cA$-bilinear with respect to the inner bimodule structure on $\cA \otimes \cA$, which then inherits a bimodule structure from the outer bimodule structure on $\cA \otimes \cA$.

There is a multi-variate version of the bimodule dual. 
For each $x,y \in \Ob(\cA)$, consider the $\cA^{\otimes n}$-bimodule ${}_{y}(\cA^{\otimes n+1})_{x} \in \Mod((\cA^{\otimes n})^e)$ with positions of $\cA$-multiplications specified as follows:
\begin{equation}  \label{A_multi_outer}
{}_{y}(\cA^{\otimes n+1})_{x} \, = \, {}_{y}\cA_{R_1} \otimes {}_{L_1}\cA_{R_2} \otimes \ldots \otimes {}_{L_{n-1}}\cA_{R_n} \otimes {}_{L_n} \cA_{x}
\end{equation}

\bdf
Given bimodules $M_1,\ldots,M_n \in \Mod(\cAe)$, then we define their \emph{bimodule multidual} to be the bimodule $\MD(M_1,\ldots,M_n) \in \Mod(\cAe)$ defined by
\begin{equation*}
	{}_{y}(\MD(M_1,\ldots,M_n))_{x} \, := \, \Homcom_{(\cA^{\otimes n})^e}( \, {}_{L_1} (M_1)_{R_1} \otimes \ldots \otimes {}_{L_n} (M_n)_{R_n} \, , \, {}_{y}(\cA^{\otimes n+1})_{x} \, )
\end{equation*}
\edf

We will write $(\cA^{\otimes n+1})_{{\rm outer}}$ to emphasize the ``outer bimodule structure'' of \eqref{A_multi_outer} as one varies $x$ and $y$. Then it is easy to see that there is an isomorphism
\begin{equation}  \label{MD_star_pre}
(\cA^{\otimes n+1})_{{\rm outer}} \otimes_{\cA} (\cA^{\otimes m+1})_{{\rm outer}} \, \cong \, (\cA^{\otimes m+n+1})_{{\rm outer}}
\end{equation}
which also preserves the ``inner $\cA^{\otimes m+n}$-bimodule structure''.

As a consequence, for any bimodules $M_1,\ldots,M_{m+n} \in \Mod(\cAe)$, there is a canonical map of $\cA$-bimodules:
\begin{equation}   \label{MD_star}
	\ast \, : \, \MD(M_1,\ldots,M_n) \otimes_{\cA}  \MD(M_{n+1},\ldots,M_{n+m}) \raq \MD(M_1,\ldots,M_{m+n})
\end{equation}

Moreover, this map is associative: $
	( \alpha \ast \beta ) \ast \gamma  =  \alpha \ast (\beta \ast \gamma)$
for $\alpha \in {}_z(\MD(M_1,\ldots,M_n))_y$, $\beta \in {}_y(\MD(M_{n+1},\ldots,M_{n+m}))_x$,  $\gamma \in {}_x(\MD(M_{n+m+1},\ldots,M_{n+m+l}))_w$. 

Denote by ${}_{\tau}(\cA^{\otimes n})_{\id}$ the $\cA^{\otimes n}$-bimodule where the left action is twisted by the cyclic rotation map $\tau : \cA^{\otimes n} \xra{\cong} \cA^{\otimes n}$. More precisely, the positions of $\cA$-multiplications are specified as follows:
\begin{equation*}  
	{}_{\tau}(\cA^{\otimes n})_{\id} \, = \, {}_{L_n}\cA_{R_1} \otimes {}_{L_1}\cA_{R_2} \otimes \ldots \otimes {}_{L_{n-1}}\cA_{R_n} 
\end{equation*}

\bdf
Given bimodules $M_1,\ldots,M_n \in \Mod(\cAe)$, then we define their \emph{naturalized multidual} to be the cochain complex 
$\MD_{\natural}(M_1,\ldots,M_n) \in \Ch(k)$ defined by
\begin{equation}  \label{MD_nat_eq}
	\MD_{\natural}(M_1,\ldots,M_n) \, := \, \Homcom_{(\cA^{\otimes n})^e}( \, {}_{L_1} (M_1)_{R_1} \otimes \ldots \otimes {}_{L_n} (M_n)_{R_n} \, , \, {}_{\tau}(\cA^{\otimes n})_{\id} \, )
\end{equation}
\edf

Notice that there is an isomorphism of $\cA^{\otimes n}$-bimodules
\begin{equation}   \label{MD_nat_pre}
	((\cA^{\otimes n+1})_{{\rm outer}})_{\natural} \, \cong \, {}_{\tau}(\cA^{\otimes n})_{\id}
\end{equation}

As a consequence, there is a canonical map of cochain complexes
\begin{equation}   \label{MD_nat}
	(\MD(M_1,\ldots,M_n))_{\natural} \raq \MD_{\natural}(M_1,\ldots,M_n)
\end{equation}

\blm  \label{MD_map_isom_1}
If $M_1,\ldots,M_{m+n}$ are projective of finite rank as graded $\cA$-bimodules, then the map \eqref{MD_star} is an isomorphism. 
If $M_1,\ldots,M_{n}$ are projective of finite rank as graded $\cA$-bimodules, then the map \eqref{MD_nat} is an isomorphism. 
\elm

\bpf
Recall that \eqref{MD_star_pre} leads to \eqref{MD_star}, and \eqref{MD_nat_pre} leads to \eqref{MD_nat}, both by putting the tensor operation inside the codomain position of $\Homcom$. This gives an isomorphism if the domain position consists of objects that are projective of finite rank.
\epf

In the case when $M_1 = \ldots = M_n = M$, we will denote by
\begin{equation*}
	\begin{split}
		\MD^{(n)}(M) \, &:= \, \MD(M, \stackrel{(n)}{\ldots}, M) \, = \, 
		\Homcom_{(\cA^{\otimes n})^e}(M^{\otimes n} , {}_{\bullet}(\cA^{\otimes n+1})_{\bullet}) \\ 
		\MD^{(n)}_{\natural}(M) \, &:= \, \MD_{\natural}(M, \stackrel{(n)}{\ldots}, M) \, = \, 
		\Homcom_{(\cA^{\otimes n})^e}(M^{\otimes n} , {}_{\tau}(\cA^{\otimes n})_{\id}) 
	\end{split}
\end{equation*}

Notice that $\MD^{(n)}_{\natural}(M)$ has an action by the cyclic group $C_n$, where the generator $\tau$ acts on the Hom complex $\Homcom_{(\cA^{\otimes n})^e}(M^{\otimes n} , {}_{\tau}(\cA^{\otimes n})_{\id}) $ by $\tau(f) = \tau \circ f \circ \tau^{-1}$.

By \eqref{MD_star}, \eqref{MD_nat} and Lemma \ref{MD_map_isom_1}, we have
\bcor  \label{tensor_of_dual_cor}
There is a canonical $C_n$-equivariant map of cochain complexes
\begin{equation*} 
	(M^{\vee} \otimes_{\cA} \stackrel{(n)}{\ldots} \otimes_{\cA} M^{\vee})_{\natural} \raq \MD^{(n)}_{\natural}(M)
\end{equation*}
which is an isomorphism if $M$ is projective of finite rank as a graded $\cA$-bimodule.
\ecor

\subsection{The higher Hochschild complexes}

\bdf
Given a $k$-flat small dg category $\cA$, its \emph{$r$-th higher cohomological Hochschild complex} ($r \geq 1$) is the complex
\begin{equation}  \label{CH_r_A}
	C_H^{(r)}(\cA) \, := \, \MD^{(r)}_{\natural}(\cR(\cA))
\end{equation}
\edf

We now describe it in more explicit terms. Again denote by $\scO = \Ob(\cA)$.
For any tuple $\vec{x} = (x_0,\ldots,x_p)$ of objects $x_j \in \scO$, we write
\begin{equation*}
	\cA(\vec{x}_i) \, := \, 
	\begin{cases*}
		{}_{x_p}\cA_{x_{p-1}} \otimes \ldots \otimes {}_{x_1}\cA_{x_0} & if $p > 0$ \\
		k & if $p = 0$
	\end{cases*}
\end{equation*}

For each $p_1,\ldots,p_r \geq 0$, and for any tuple $\mathbf{x} = (\vec{x}_1, \ldots, \vec{x}_r)$, where $\vec{x}_i = (x_{i,0},\ldots,x_{i,p_i})$ for $x_{i,j} \in \scO$, we write
\begin{equation}  \label{HHH1}
	\begin{split}
		C_H^{[p_1,\ldots,p_r]}(\cA)_{\mathbf{x}} \, &:= \, \Homcom_{k}(\, \cA(\vec{x}_1)  \otimes \ldots \otimes \cA(\vec{x}_r) \, , \, {}_{x_{r,p_r}}\cA_{x_{1,0}} \otimes \ldots \otimes {}_{x_{r-1,p_{r-1}}}\cA_{x_{r,0}} \, )	\\
		C_H^{[p_1,\ldots,p_r]}(\cA) \, &:= \, \prod_{\mathbf{x} = (\vec{x}_1, \ldots, \vec{x}_r)} C_H^{[p_1,\ldots,p_r]}(\cA)_{\mathbf{x}}
	\end{split}
\end{equation} 

Then one can see that (recall the notation from \eqref{RnA_def})
\begin{equation*}
		C_H^{[p_1,\ldots,p_r]}(\cA) \, = \, \Homcom_{(\cA^{\otimes n})^e}(\, \cR_{p_1}(\cA) \otimes \ldots \otimes \cR_{p_r}(\cA) \, , \, {}_{\tau}(\cA^{\otimes r})_{\id} \, ) 
\end{equation*} 
so that $C_H^{[\bullet,\ldots,\bullet]}(\cA)$ forms an $r$-fold cosimplicial system $\Delta \times \stackrel{(r)}{\ldots} \times \Delta \ra \Ch(k)$, because each $\cR_{\bullet}(\cA)$ forms a simplicial system $\Delta^{\op} \ra \Mod(\cAe)$. 
The higher Hochschild complex $C_H^{(r)}(\cA)$ can be realized as the direct product total complex of this $r$-fold cosimplicial system:
\begin{equation}   \label{HHH3}
	C_H^{(r)}(\cA) \, = \, \prod_{(p_1,\ldots,p_r) \in \bN^r} \, C_H^{[p_1,\ldots,p_r]}(\cA)[-p_1-\ldots -p_r]
\end{equation}

We will consider the $m$-shifted version of \eqref{CH_r_A}:
\begin{equation}  \label{CH_r_A_shifted}
	C_H^{(r)}(\cA,m) \, := \, \MD^{(r)}_{\natural}(\cR(\cA)[m])
\end{equation}

Recall from above that $\MD^{(r)}_{\natural}(M)$ has an action by the cyclic group $C_r$. In particular, we may take the $C_r$-invariants of $C_H^{(r)}(\cA,m)$, which will be denoted as $C_H^{(r)}(\cA,m)^{C_r}$, and called the \emph{$m$-shifted cyclic invariant $r$-th higher cohomological Hochschild complex}. (Of course, $C_H^{(r)}(\cA,m)$ is isomorphic to $C_H^{(r)}(\cA)[-mr]$, but we take the shift as in \eqref{CH_r_A_shifted} in order to specify the $C_r$-action. In particular, $C_H^{(r)}(\cA,m)^{C_r} \neq C_H^{(r)}(\cA)^{C_r}[-mr]$ in general.)

The cyclic invariant higher cohomological Hochschild complexes have an alternative interpretation. 
We now assume that $\cA$ consists of the data of a set $\scO$ together with a collection of graded $k$-modules ${}_y\cA_x = \cA(x,y)$ for each $x,y \in \scO$.
Thus, $\cA$ is not yet given the structure of a dg category.
Then \eqref{HHH1} and \eqref{HHH3} are still well-defined, but now only as a graded $k$-module (we no longer have an $n$-cosimplicial system $C_H^{[\bullet,\ldots,\bullet]}(\cA)$, but the direct product \eqref{HHH3} is still well-defined). 
The graded $k$-module $C_H^{(r)}(\cA,m) = C_H^{(r)}(\cA)[-mr]$, together with the $C_r$-action, is still well-defined. 

\bdf  \label{disk_def}
An \emph{$\scO$-colored disk} consists of the following data:
\begin{enumerate}
	\item A finite set of punctures $\Sigma \subset \partial \bD$ on the boundary of the $2$-dimensional disk $\bD$, with a partition $\Sigma = \Sigma^+ \amalg \Sigma^-$. The punctures in $\Sigma^+$ will be called an input puncture, and those in $\Sigma^-$ will be called an output puncture.
	\item A map of sets $o : \pi_0(\partial \bD \setminus \Sigma) \ra \scO$. 
\end{enumerate}

A diffeomorphism of a $\scO$-colored disk is a diffeomorphism of the disk $\bD$ presreving both data (1) and (2). Define the category $\frf$ whose objects are $\scO$-colored disks, and whose morphisms are isotopy classes of diffeomorphisms. Clearly, $\frf$ is a groupoid.

For each $\zeta \in \Sigma$, denote by $\zeta', \zeta'' \in \pi_0(\partial \bD \setminus \Sigma)$ the two connected components such that $(\zeta'', \zeta, \zeta')$ is in counterclockwise direction in a neighborhood of $\zeta$. Define $\nu(\zeta) \in \scO \times \scO$ by
\begin{equation*}
	\nu(\zeta) \, = \,
	\begin{cases*}
		(o(\zeta''),o(\zeta')) & if $\zeta \in \Sigma^+$ \\
		(o(\zeta'),o(\zeta'')) & if $\zeta \in \Sigma^-$
	\end{cases*}
\end{equation*}
\edf

Given $r \geq 1$, $p_1,\ldots,p_r \geq 0$, and a tuple $\mathbf{x} = (\vec{x}_1, \ldots, \vec{x}_r)$, where $\vec{x}_i = (x_{i,0},\ldots,x_{i,p_i})$ for $x_{i,j} \in \scO$, denote by $\bD(\mathbf{x})$ the $\scO$-colored disk with $n$ output punctures, denoted as $\zeta^-_1,\ldots, \zeta^-_r$, put in counterclockwise order, and with $p_i$ input punctures in the counterclockwise segment of $\partial \bD$ from $\zeta^-_i$ to $\zeta^-_{i+1}$ (where we write $\zeta^-_{r+1} := \zeta^-_1$). Thus, $\pi_0(\partial \bD \setminus \Sigma)$ has $p_i + 1$ elements in the counterclockwise segment of $\partial \bD$ from $\zeta^-_i$ to $\zeta^-_{i+1}$. These are to be colored by the tuple $\vec{x}_i$, again in counterclockwise direction. For instance, the following is an example for $r = 3$ and $(p_1,p_2,p_3) = (1,0,2)$:
\begin{equation*}
	\includegraphics*[scale=0.3]{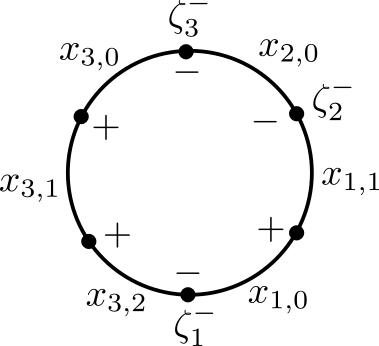}
\end{equation*}

This allows us to give a combinatorial description of $\frf$. 
Notice that the cyclic group $C_r$ acts on the set $\{ ((p_1,\ldots,p_r),\mathbf{x}) \}$ of data described above.
Define $\frf_r$ to be the action groupoid. {\it i.e.,} $\Ob(\frf_r) = \{ ((p_1,\ldots,p_r),\mathbf{x}) \}$, and $\Hom_{\frf_r}((\vec{p},\mathbf{x}),(\vec{p}',\mathbf{x}')) = \{ \sigma \in C_r | \sigma \cdot (\vec{p},\mathbf{x}) = (\vec{p}',\mathbf{x}') \}$.
Define $\frf_{\geq 1} = \coprod_{n \geq 1} \frf_n$. Then the above construction of $\bD(\mathbf{x})$ defines a functor $\frf_{\geq 1} \ra \frf$, which is fully faithful, and whose essential image consists of $\scO$-colored disks with at least one output puncture.

The complexes \eqref{HHH1} can also be described by using the language of $\scO$-colored disk. For example, we have
\begin{equation}  \label{HHH4}
	C_H^{[p_1,\ldots,p_r]}(\cA)_{\mathbf{x}} \, = \, \Endcom^{{\rm rd}}(\cA)(\bD(\mathbf{x})) \, := \,  \Homcom_k \Biggl( \, \bigotimes_{\zeta \in \Sigma^+} \, \cA(\nu(\zeta)) \, , \,  \bigotimes_{\zeta \in \Sigma^-} \, \cA(\nu(\zeta)) \, \Biggr)
\end{equation}

The right hand side of \eqref{HHH4} depends only on the geometry of $\bD(\mathbf{x})$, so that it is functorial in $\bD(\mathbf{x}) \in \frf_{\geq 1}$ (in fact it is functorial over $\frf$). More generally, we can consider its $(a,b)$-shifted version for any $a,b \in \bZ$:
\begin{equation}  \label{End_rd_ab}
	\begin{split}
		\Endcom^{{\rm rd}}_{(a,b)}(\cA) \, &: \, \frf \raq \grMod(k) \\
		\Endcom^{{\rm rd}}_{(a,b)}(\cA)( \bD(\mathbf{x}) )	\, &:= \, \Homcom_k \Biggl( \, \bigotimes_{\zeta \in \Sigma^+} \, \Bigl( \cA(\nu(\zeta))[a] \Bigr) \, , \,  \bigotimes_{\zeta \in \Sigma^-} \, \Bigl( \cA(\nu(\zeta))[-b] \Bigr) \, \Biggr)[a+b]
	\end{split}
\end{equation}

Now we will fix $m \in \bZ$, and consider the case $(a,b) = (1,m)$.
Then we have
\begin{equation}  \label{End_rd_lim}
	C_H^{(r)}(\cA,m)^{C_r}[m+1] \, = \, \lim \, \Bigl(\, \Endcom^{{\rm rd}}_{(1,m)}(\cA)|_{\frf_r} \, : \, \frf_r \raq \grMod(k) \, \Bigr)
\end{equation}

Define the $k$-linear map
\begin{equation*}  
	\circ \, : \, C_H^{(p)}(\cA,m)^{C_p}[m+1] \, \otimes \, C_H^{(q)}(\cA,m)^{C_q}[m+1] \raq C_H^{(p+q-1)}(\cA,m)^{C_{p+q-1}}[m+1]
\end{equation*}
by the formula (see explanation below):
\begin{equation}  \label{F_circ_G_1}
	\includegraphics*[scale=0.2]{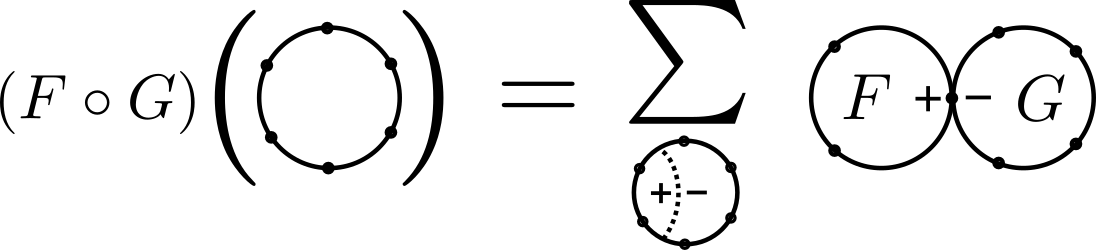}
\end{equation}

Given $\scO$-colored disks $\bD(\mathbf{x}')$ and $\bD(\mathbf{x}'')$, if we choose a pair $\zeta^+ \in \Sigma^+( \bD(\mathbf{x}') )$ and $\zeta^- \in \Sigma^-( \bD(\mathbf{x}'') )$ such that $\nu(\zeta^+) = \nu(\zeta^-)$, then one may take a connected sum of $\bD(\mathbf{x}')$ and $\bD(\mathbf{x}'')$ along $\zeta^+$ and $\zeta^-$. Denote the resulting $\scO$-colored disk by $\bD = \bD(\mathbf{x}') {}_{\zeta^+}\#_{\zeta^-} \bD(\mathbf{x}'')$. 
Then there is an obvious composition map
\begin{equation}  \label{End_rd_comp}
	{}_{\zeta^+}\circ_{\zeta^-} \, : \, \Endcom^{{\rm rd}}_{(1,m)}(\cA)(\bD(\mathbf{x}')) \, \otimes \, \Endcom^{{\rm rd}}_{(1,m)}(\cA)(\bD(\mathbf{x}'')) \raq \Endcom^{{\rm rd}}_{(1,m)}(\cA)(\bD))
\end{equation}
defined by putting the output of $G \in \Endcom^{{\rm rd}}_{(1,m)}(\cA)(\bD(\mathbf{x}''))$ at the $\zeta^-$ position into the input of $F \in \Endcom^{{\rm rd}}_{(1,m)}(\cA)(\bD(\mathbf{x}'))$ at the $\zeta^+$ position. (See \cite{Yeu5} for the precise Koszul signs involved.)

Now suppose we are given $F \in C_H^{(p)}(\cA,m)^{C_p}$ and $G \in C_H^{(q)}(\cA,m)^{C_q}$, then we define $F \circ G \in C_H^{(p+q-1)}(\cA,m)^{C_{p+q-1}}$ by
\begin{equation}  \label{F_circ_G_2}
	(F \circ G)(\bD) \, = \, \sum_{\bD = \bD(\mathbf{x}') {}_{\zeta^+}\#_{\zeta^-} \bD(\mathbf{x}'')} \, F( \bD(\mathbf{x}') ) \,\, {}_{\zeta^+} \! \circ_{\zeta^-} G( \bD(\mathbf{x}'') )
\end{equation}
where we sum over all the possible ways (up to diffeomorphism) to write $\bD$ as a connected sum $\bD(\mathbf{x}') {}_{\zeta^+}\#_{\zeta^-} \bD(\mathbf{x}'')$ of two $\scO$-colored disks. (Since $F$ and $G$ are homogeneous of weight $p$ and $q$ respectively, we regard the term in the sum \eqref{F_circ_G_2} to be zero unless $\bD(\mathbf{x}')$ has $p$ output punctures and $\bD(\mathbf{x}'')$ has $q$ output punctures.)

Notice that, given $\bD$, the different ways (up to diffeomorphism) of writing $\bD$ as a connected sum $\bD = \bD(\mathbf{x}') {}_{\zeta^+}\#_{\zeta^-} \bD(\mathbf{x}'')$ is in a natural bijection to $\pi_0(\partial \bD \setminus \Sigma) \times \pi_0(\partial \bD \setminus \Sigma)$. Namely, given $b',b'' \in \pi_0(\partial \bD \setminus \Sigma)$, draw a line from $b'$ to $b''$, and write ``$-$'' on the right hand side, and ``$+$'' on the left hand side, as in the summation index of \eqref{F_circ_G_1}. This splits $\bD$ as a connected sum along this line. 
Thus, \eqref{F_circ_G_1} is a pictorial description of \eqref{F_circ_G_2}.

The following result appeared in \cite{KV13, KTV}, and independently discovered by us. We give a proof that neglects signs. A more careful proof can be found in \cite{Yeu5}:
\bthm  \label{Lie_bracket_thm}
Define the bracket 
\begin{equation*}
	\{-,-\} \, : \, C_H^{(p)}(\cA,m)^{C_p}[m+1] \, \otimes \, C_H^{(q)}(\cA,m)^{C_q}[m+1] \raq C_H^{(p+q-1)}(\cA,m)^{C_{p+q-1}}[m+1]
\end{equation*}
by the formula
\begin{equation*}
	\{ F,G \} \, := \, F \circ G - (-1)^{|F||G|} G \circ F
\end{equation*}
then $\bigoplus_{r \geq 1} C_H^{(r)}(\cA,m)^{C_r}[m+1]$ becomes a graded Lie algebra.
\ethm

\bpf[Sketch of proof (see \cite{Yeu5} for details).]
Evaluate the sum 
\begin{equation}  \label{FGH_cyclic_sum}
	\{ F, \{G,H\} \} + (-1)^{|F|(|G|+|H|)} \{ G, \{H,F\} \} + (-1)^{|H|(|F|+|G|)} \{ H, \{F,G\} \}
\end{equation}
on an $\scO$-colored disk $D$. The terms in the resulting sum is indexed by the data of (i) a splitting of the disk by two non-intersecting dotted line, as in the summation index of \eqref{F_circ_G_1}; (ii) for each dotted line, a choice of which side to designate $+$ (and $-$ to the other side); (iii) a choice of filling in $F,G,H$ to the resulting three disks. Each such terms appears twice in \eqref{FGH_cyclic_sum}, cancelling each other. For example, the term corresponding to
\begin{equation*}
	\includegraphics*[scale=0.25]{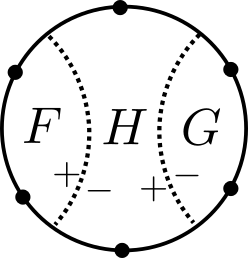}
\end{equation*}
appears in both $\{ F, \{G,H\}\}$ and $(-1)^{|F|(|G|+|H|)} \{ G, \{H,F\} \}$.
\epf

\brm
In \cite{Yeu5}, we defined a \emph{graded pseudo pre-Lie algebra} to be a graded $k$-module $L$ together with maps 
\begin{equation*}
	\begin{split}
		\circ \, &: \, L \otimes L \raq L \\
		m_{2\ra 1 \la 3} \, &: \, L \otimes L \otimes L \raq L \\
		m_{1\la 3 \ra 2} \, &: \, L \otimes L \otimes L \raq L
	\end{split}
\end{equation*}
such that
\begin{equation*}
	\begin{split}
		m_{2\ra 1 \la 3}(x,y,z) \, &= \, (-1)^{|y||z|} m_{2\ra 1 \la 3}(x,z,y) \\
		m_{1\la 3 \ra 2}(x,y,z) \, &= \, (-1)^{|x||y|} m_{1\la 3 \ra 2}(y,x,z) \\
		(x \circ y) \circ z - x \circ (y\circ z) \, &= \, m_{2\ra 1 \la 3}(x,y,z) - m_{1\la 3 \ra 2}(x,y,z)
	\end{split}
\end{equation*}

Notice that a graded pre-Lie algebra is precisely a graded pseudo pre-Lie algebra with $m_{1\la 3 \ra 2} = 0$. It can be shown that if $L$ is a graded pseudo pre-Lie algebra, then the bracket $\{F,G\} = F \circ G - (-1)^{|F||G|} G \circ F$ gives a graded Lie algebra structure on $L$.

In \cite{Yeu5}, we defined a graded $\scO$-colored ribbon dioperad to be a functor $\cP : \frf \ra \grMod(k)$ (an example is \eqref{End_rd_ab}), together with composition maps as in \eqref{End_rd_comp}, that satisfies certain associativity conditions. It is proved in \cite{Yeu5} that for any graded $\scO$-colored ribbon dioperad, the formula \eqref{F_circ_G_2} is part of the structure of a graded pseudo pre-Lie algebra on the limit total object $\bigoplus_{n \geq 0} \lim(\cP : \frf_n \ra \grMod(k))$.
\erm

Recall that, so far, we have assumed that $\cA$ consists only of the data of a set $\scO$ together with a collection of graded $k$-modules ${}_y\cA_x = \cA(x,y)$ for each $x,y \in \scO$.
We are now ready to consider dg category (in fact $A_{\infty}$-category) structures on $\cA$.

Notice that the bracket in Theorem \ref{Lie_bracket_thm} has weight grading $-1$. Thus, in the graded Lie algebra $\bigoplus_{r \geq 1} C_H^{(r)}(\cA,m)^{C_r}[m+1]$, the part $r = 1$ is a Lie subalgebra, and is independent of $m$. This is precisely the underlying graded $k$-module of the usual Hochschild cochain complex 
\begin{equation*}
C_H(\cA)[1] \, = \, \prod_{p \geq 0} \,\,\, \prod_{\vec{x} = (x_0,\ldots,x_p) \in \scO^{p+1}} \Homcom_k( \, \cA(\vec{x})[p] \, , \, {}_{x_p}\cA_{x_0}[1] \, )
\end{equation*}

It is clear that the bracket $\{-,-\}$ in Theorem \ref{Lie_bracket_thm}, when restricted to this weight $1$ part, is precisely the usual Gerstenhaber bracket. Recall that a (possibly curved, and not necessarily unital) $A_{\infty}$-category structure is precisely a Maurer-Cartan element $\mu = \mu_0 + \mu_1 + \mu_2 + \ldots$ on this graded Lie algebra $(C_H(\cA)[1],\{-,-\})$. 
We will only consider the cases where $(\cA,\mu)$ is non-curved, {\it i.e.,} $\mu_0 = 0$.

From now on, we assume that $(\cA,\mu)$ is a non-curved (not necessarily unital) $A_{\infty}$-category. Since $C_H(\cA)[1]$ is a graded Lie subalgebra of $\bigoplus_{r \geq 1} C_H^{(r)}(\cA,m)^{C_r}[m+1]$, the element $\mu$ is still a Maurer-Cartan element in the entire graded Lie algebra $\bigoplus_{r \geq 1} C_H^{(r)}(\cA,m)^{C_r}[m+1]$. Thus, we can use $\mu$ to define a differential $d = d_{\mu} = \{\mu,-\}$, so that there is a dg Lie algebra
\begin{equation}  \label{cyclic_inv_HHH_dgla}
 \Bigl( \, \bigoplus_{r \geq 1} C_H^{(r)}(\cA,m)^{C_r}[m+1] \, , \, \{-,-\}   \, , \,  d_{\mu} \, \Bigr)
\end{equation}
Again, since the bracket $\{-,-\}$ has weight grading $-1$, this differential $d_{\mu}$ preserves the weight grading, so that each weight component of \eqref{cyclic_inv_HHH_dgla} is a subcomplex. We call the complex $(C_H^{(r)}(\cA,m)^{C_r}, d_{\mu})$ the \emph{$m$-shifted cyclic invariant $r$-th higher cohomological Hochschild complex} of the $A_{\infty}$-category $(\cA,\mu)$. If $(\cA,\mu)$ is in fact a unital dg category, then it coincides with the $C_r$-invariants of \eqref{CH_r_A_shifted}.

\brm
In our above discussion, we have avoided the case $r = 0$, but all our discussions extend over this case. For example, the analogue of $C_H^{(r)}(\cA)^{C_r}$ for $r=0$ is the $k$-linear dual of the Connes complex $C^{\lambda}(\cA)$ (see, {\it e.g.}, \cite{Lod98} or \cite{Yeu3} for the notation $C^{\lambda}(\cA)$).
\erm

Since the bracket of \eqref{cyclic_inv_HHH_dgla} has weight $-1$, the weight $\geq 2$ part is a dg Lie ideal, and hence in particular a dg Lie subalgebra. We will consider its completion
\begin{equation}  \label{cyclic_inv_HHH_dgla_2}
	\Bigl( \, \prod_{r \geq 2} C_H^{(r)}(\cA,m)^{C_r}[m+1] \, , \, \{-,-\}   \, , \,  d_{\mu} \, \Bigr)
\end{equation}

\bdf  \label{preCY_def_1}
Assume that $k$ is a field of characteristic $0$, and $(\cA,\mu)$ is a non-curved (not necessarily unital) $A_{\infty}$-category. Then an \emph{$n$-pre-Calabi-Yau structure} is a Maurer-Cartan element in the dg Lie algebra \eqref{cyclic_inv_HHH_dgla_2} for $m = 2-n$.
In other words, it consists of $\pi = \pi_2 + \pi_3 + \ldots$, where $\pi_r \in C_H^{(r)}(\cA,m)^{C_r}[m+1]$ is a degree $1$ element, such that $d_{\mu}(\pi) + \frac{1}{2}\{\pi,\pi \} = 0$.
\edf

If we write $\pi_1 = \mu$, and $\widetilde{\pi} = \pi_1 + \pi_2 + \pi_3 + \ldots$, then the Maurer-Cartan element $d_{\mu}(\pi) + \frac{1}{2}\{\pi,\pi \} = 0$ simplifies to $\{ \widetilde{\pi} , \widetilde{\pi} \} = 0$. Thus, the $A_{\infty}$-category structure and the pre-Calabi-Yau structure can be combined into the following

\bdf  \label{preCY_def_2}
Assume that $k$ is a field of characteristic $0$, then 
a (possibly curved, not necessarily unital) \emph{$n$-pre-Calabi-Yau category} over $k$ consists of a set $\scO$ with a collection of graded $k$-vector spaces ${}_y\cA_x = \cA(x,y)$ for each $x,y \in \scO$, together with a Maurer-Cartan element in the graded Lie algebra
\begin{equation*}  
	\Bigl( \, \prod_{r \geq 1} C_H^{(r)}(\cA,m)^{C_r}[m+1] \, , \, \{-,-\}   \, \Bigr)
\end{equation*}
\edf

\brm
In Definitions \ref{preCY_def_1} and \ref{preCY_def_2}, we have assumed that $k$ is a field of characteristic $0$. The assumption that $k$ is a field is imposed to guarantee that each cochain complex $(\cA(x,y),d)$ is automatically cofibrant over $k$, which ensures that Hochschild (co)homology have the expected properties. The assumption that $k$ has characteristic $0$ ensures that the procedure of taking $C_r$-invariants in $C_H^{(r)}(\cA,m)$ is well-behaved. In general, there are a priori $4$ possible versions: one either take invariants, coinvariants, homotopy invariants, or homotopy coinvariants, all of which coincide in characteristic $0$. See \cite{Yeu5} for some discussions along this line.
\erm


\subsection{The extended necklace dg Lie algebra}  \label{sec_ext_necklace}

Now we return to the usual setting, where $(\cA,\mu)$ is a unital dg category over a commutative unital ring $k$. In the beginning of Section \ref{sec_nc_DR}, we used the principles \eqref{NC_prin_cSA} and \eqref{NC_prin_naturalize} to justify \eqref{scX_n_def} as a noncommutative analogue of global $n$-forms. By a completely analogous reasoning, one may regard the following as a noncommutative analogue of global $r$-polyvectors:
\begin{equation*}
	 (\cSA^{\vee} \otimes_{\cA} \stackrel{(r)}{\ldots} \otimes_{\cA} \cSA^{\vee})_{\natural, C_r}
\end{equation*}

However, we will make $3$ modifications to this:
\begin{enumerate}
	\item By Lemma \ref{tensor_of_dual_cor}, the complex  
	$(\cSA^{\vee} \otimes_{\cA} \stackrel{(r)}{\ldots} \otimes_{\cA} \cSA^{\vee})_{\natural}$ is $C_r$-equivariantly isomorphic to $\MD^{(r)}_{\natural}(\cSA)$ if $\cSA$ is projective of finite rank as a graded bimodule. In general, when the two differ, we regard $\MD^{(r)}_{\natural}(\cSA)$ to be the ``correct'' version to consider, as there are pathologies involved in taking the tensor products of duals.
	\item We will take $C_r$-invariants instead of $C_r$-coinvariants. They coincide in characteristic $0$.
	\item We will take an $m$-shifted version.
\end{enumerate}

Combining these considerations, we take
\begin{equation*}
	\scP^{(r)}(\cA;m) \, := \, \MD^{(r)}_{\natural}(\cSA[m])^{C_r}
\end{equation*}

We compare this with the $m$-shifted cyclic invariant higher cohomological Hochschild complexes
\begin{equation*}
	C_H^{(r)}(\cA,m)^{C_r} \, = \, \MD^{(r)}_{\natural}(\cR(\cA)[m])^{C_r}
\end{equation*} 

The canonical surjective quasi-isomorphism $\cR(\cA) \xronto{\sim} \cSA$ induces an injective map of complexes
\begin{equation}  \label{ext_necklace_to_HHH}
	\bigoplus_{r \geq 1} \, \scP^{(r)}(\cA;m)[m+1] \rintoq \bigoplus_{r \geq 1} \, C_H^{(r)}(\cA,m)^{C_r}[m+1] 
\end{equation}

\bthm  \label{ext_necklace_Lie_subalg}
The subcomplex \eqref{ext_necklace_to_HHH} is a dg Lie subalgebra. 
Moreover, if $\cA$ is almost cofibrant and $\bQ \subset k$, then the inclusion map \eqref{ext_necklace_to_HHH} is a quasi-isomorphism.
\ethm

\bpf
Recall from \eqref{HHH3} that an element $F \in C_H^{(r)}(\cA,m)^{C_r}$ consists of a collection of elements $F_{(\vec{p},\mathbf{x})} \in C_H^{[p_1,\ldots,p_n]}(\cA)_{\mathbf{x}}$, one for each $(\vec{p},\mathbf{x})$, that is $C_r$-invariant with respect to a suitable Koszul rule.
Since $(\cA,\mu)$ is a dg category, the differential $d = d_{\mu}$ on $C_H^{(r)}(\cA,m)^{C_r}$ decomposes as $d = d_1 + d_2$, where $d_1$ is induced by the differential in $\cA$, and $d_2$ is induced by the composition maps in $\cA$.

An element $F \in C_H^{(r)}(\cA,m)^{C_r}$ is in the subcomplex \eqref{ext_necklace_to_HHH} if and only if
\begin{enumerate}
	\item $F_{(\vec{p},\mathbf{x})} = 0$ whenever $p_i > 1$ for some $1 \leq i \leq r$.
	\item $(d_2 F)_{(\vec{p},\mathbf{x})} = 0$ whenever $p_i > 1$ for some $1 \leq i \leq r$.
\end{enumerate}
Indeed, notice that assuming (1), then it suffices to check condition (2) for the case when $p_i = 2$ and $p_j \leq 1$ for $j \neq i$. In this case, the condition  $(d_2 F)_{(\vec{p},\mathbf{x})} = 0$ is a precise way of saying that the $i$-th input of $F$ satisfies a derivation property.

Assume that $F,G$ both satisfy (1),(2) above. It is clear that $\{F,G\}$ also satisfies (1). To show that it also satisfies (2), notice that $\bigoplus_{r \geq 1} C_H^{(r)}(\cA,m)^{C_r}[m+1]$ is still a dg Lie algebra with $d$ replaced by $d_2$. (This is because the graded category $(\cA,\mu_2)$ is also an $A_{\infty}$-category.) Thus, we have
$d_2(\{F,G\}) = \{ d_2(F),G \} + (-1)^{|F|} \{F,d_2(G) \}$. All the four terms $F,G,d_2(F),d_2(G)$ inside the brackets vanishes on $(\vec{p},\mathbf{x})$ whenever $p_i > 1$ for some $1 \leq i \leq r$, hence the same is true for their brackets. This shows that $\{F,G\}$ satisfies (2).
This completes the proof that \eqref{ext_necklace_to_HHH} is a dg Lie subalgebra.

For the second statement, notice that if $\cA$ is $k$-flat, then the functor $\Mod(\cAe) \ra \Ch(k)$, $M \mapsto \MD^{(r)}_{\natural}(M)$ sends quasi-isomorphisms between cofibrant bimodules to quasi-isomorphisms. If $\cA$ is almost cofibrant, then $\cR(\cA)[m] \xronto{\sim} \cSA[m]$ is a quasi-isomorphism between cofibrant bimodules, so that $\MD^{(r)}_{\natural}(\cSA[m]) \rinto \MD^{(r)}_{\natural}(\cR(\cA)[m])$ is a quasi-isomorphism. If $\bQ \subset k$, then taking $C_r$-invariants preserves quasi-isomorphisms.
\epf

\bdf
The dg Lie algebra $\bigoplus_{r \geq 1} \, \scP^{(r)}(\cA;m)[m+1]$ is called the \emph{$m$-shifted extended necklace dg Lie algebra}.
\edf

To explain the terminology, we now consider a ``framed version'' of $\scP^{(r)}(\cA;m)$, which will recover the (non-extended) necklace Lie algebra \cite{BL02} of an associative algebra.

For any subset $\scO' \subset \scO$, let $\cSAfr \in \Mod(\cAe)$ be the quotient of $\cSA$ by the sub-bimodule $\cA \otimes_{\scO'} \cA \subset \cA \otimes_{\scO} \cA \subset \cSA$. In particular, if $\scO' = \emptyset$ then $\cSAfr = \cSA$; if $\scO' = \scO$ then $\cSAfr = \Omega^1(\cA)[1]$. Let 
\begin{equation*}
		\scP^{(r)}_{\Opfr}(\cA;m) \, := \, \MD^{(r)}_{\natural}(\cSAfr[m])^{C_r}
\end{equation*}

The surjection $\cSA \ronto \cSAfr$ induces an injective map of complexes
\begin{equation}  \label{necklace_to_ext_necklace}
	\bigoplus_{r \geq 1} \, \scP^{(r)}_{\Opfr}(\cA;m)[m+1] \rintoq 	\bigoplus_{r \geq 1} \, \scP^{(r)}(\cA;m)[m+1]
\end{equation}

\bpp  \label{necklace_Lie_subalg}
The subcomplex \eqref{necklace_to_ext_necklace} is a dg Lie subalgebra.
\epp

\bpf
We continue to use the notation in the proof of Theorem \ref{ext_necklace_Lie_subalg}. An element $F \in C_H^{(r)}(\cA,m)^{C_r}$ is in the subcomplex $\scP^{(r)}_{\Opfr}(\cA;m)$ if and only if it satisfies (1), (2), and
\begin{enumerate}
	\item[(3)]  $F_{(\vec{p},\mathbf{x})} = 0$ whenever $p_i = 0$ and $\vec{x}_i = (x_{i,0})$ for $x_{i,0} \in \scO'$, for some $1 \leq i \leq r$.
\end{enumerate}
This condition is clearly preserved under $\{-,-\}$.
\epf

\bdf
In the case $\scO'=\scO$, the complex $\scP^{(r)}_{\Opfr}(\cA;m)$ will be simply denoted by $\scP^{(r)}_{\fr}(\cA;m)$. 
The dg Lie algebra $\bigoplus_{r \geq 1} \, \scP^{(r)}_{\fr}(\cA;m)[m+1]$ is called the \emph{$m$-shifted necklace dg Lie algebra}.
\edf

When the differential of $\cA$ is zero, for example when $\cA = A$ is an associative algebra, then the differential of the $m$-shifted necklace dg Lie algebra is also zero. Thus, it is a graded Lie algebra in this case, and will be simply called the \emph{$m$-shifted necklace Lie algebra}. In particular, for associative algebras, the $0$-shifted version coincides with the necklace Lie algebra defined in, {\it e.g.}, \cite{BL02}.

\bcor
On an associative algebra $A$, any double Poisson bracket \cite{VdB08a} determines a $2$-pre-Calabi-Yau structure.
\ecor

\bpf
For an associative algebra, $\scP^{(r)}_{\fr}(A;0)$ is concentrated in degree $r$. Hence, a Maurer-Cartan element in $\prod_{r \geq 2} \scP^{(r)}_{\fr}(A;0)[1]$ is precisely a double Poisson bracket. Combining Theorem \ref{ext_necklace_Lie_subalg} and Proposition \ref{necklace_Lie_subalg}, the (completed) necklace Lie algebra sits as a dg Lie subalgebra
\begin{equation*}
	\prod_{r \geq 2} \scP^{(r)}_{\fr}(A;0)[1] \rintoq \prod_{r \geq 2} \scP^{(r)}(A;0)[1] \rintoq \prod_{r \geq 2} C_H^{(r)}(\cA)^{C_r}[1]
\end{equation*}
so that a Maurer-Cartan element in the former determines one in the latter.
\epf

\brm
This result appeared in \cite{IKV21}.
More generally, if $A$ is a dg algebra, then our proof also establishes the main result of \cite{FH21}.
\erm

\brm
In the beginning of this subsection, we used the principles \eqref{NC_prin_cSA} and \eqref{NC_prin_naturalize} to justify $\scP^{(r)}(\cA;m)$ as a noncommutative analogue of $m$-shifted global $p$-polyvector fields. One can also use the same principles to obtain the Lie bracket on $\bigoplus_{r \geq 1} \, \scP^{(r)}(\cA;m)[m+1]$. 

On a smooth variety $X$, polyvector fields may be regarded as global functions on the $1$-shifted cotangent bundle. The Schouten-Nijenhuis bracket on polyvector fields may be regarded as the Poisson bracket associated to the canonical $1$-shifted symplectic structure on the $1$-shifted cotangent bundle. Both of these ingredients have noncommutative analogue along the lines of Section \ref{sec_FNAG}. In particular, as a noncommutative analogue of the statement that the (shifted) cotangent bundle has a canonical symplectic structure, we have shown in \cite{Yeu1} that the $(m+1)$-Calabi-Yau completion of a smooth dg category $\cA$ has an $(m+1)$-Calabi-Yau structure, which moreover admits an explicit formula if $\cA$ is finitely cellular. Thus, at least when $\cA$ is finitely cellular, one can use this Calabi-Yau structure to construct a bracket on $\bigoplus_{r \geq 1} \, \scP^{(r)}(\cA;m)[m+1]$, which coincides with the one obtained via Theorem \ref{ext_necklace_Lie_subalg}. 

This was indeed the approach followed by the earlier version of this paper. We started with the formula for the bracket on $\bigoplus_{r \geq 1} \, \scP^{(r)}(\cA;m)[m+1]$ in the case when $\cA$ is finitely cellular, and guessed the formula in general, and then extend that to \eqref{cyclic_inv_HHH_dgla}. In this present version, we disregard the chronological order of our discovery of results, and start our presentation with Theorem \ref{Lie_bracket_thm} instead. This leads to a much cleaner presentation.
\erm

\section{Derived moduli spaces of representations}  \label{sec_DRep}

From now on, we assume that $k$ is a field of characteristic $0$. 
A $k$-category is a category enriched over $(\Vect_k, \otimes)$. 
In this section, we will study the (derived) moduli spaces of representations of small $k$-categories and more generally small dg categories. We give an example here. The conceptual underpinnings of this example will be provided in Sections \ref{Rep_subsec} and \ref{DRep_subsec} below.

A small $k$-category can be thought of as an associative algebra with many objects. 
Just like associative algebras, $k$-categories can also be defined by generators and relations. A typical example is the followng:
\begin{equation}  \label{A_V_ex1}
	A \, = \, k
	\left \langle
	\begin{tikzcd}
		x \ar[rr, bend left = 20, "f"]  \ar[loop left, distance=3em, "h"]
		& &  y \ar[ll, bend left = 20, "g"]
	\end{tikzcd} 
	\right \rangle
	\Bigg /
	\Bigl(  gf = h^2 \Bigr)
\end{equation}
where objects of $A$ are the vertices $\{x,y\}$ in the displayed quiver, and morphisms of $A$ are sums of monomials of composable arrows, modulo the ideal generated by the relation $gf-h^2$.

A representation of $A$ is a $k$-linear functor $A \ra \Vect_k^{{\rm fd}}$.
To study the space of representations of a small $k$-category $A$,
a standard approach is to choose a dimension vector 
$\vec{n} : \mathscr{O} \ra \bN$, where $\mathscr{O} = \Ob(A)$, and consider $k$-linear functors $T : A \ra \Vect_k$,
such that $T(x) = k^{\vec{n}(x)}$ for each $x \in \mathscr{O}$.
If $A$ is presented as $A = k\langle Q \rangle/I$, then to specify $T$, one should specify a matrix $T(f)$ for each $f \in Q$, such that they satisfy the matrix relations specified in $I$.
For example, for $A$ defined in \eqref{A_V_ex1}, specifying a representation of $A$ is equivalent to giving
$\bigl(  f_{ij}  \bigr) \in \bM_{\vec{n}(y) \times \vec{n}(x)}(k)$,
$\bigl(  g_{ij}  \bigr) \in \bM_{\vec{n}(x) \times \vec{n}(y)}(k)$,
$\bigl(  h_{ij}  \bigr) \in \bM_{\vec{n}(x) \times \vec{n}(x)}(k)$,
so that the equation 
$\sum_{1 \leq j \leq \vec{n}(y)} g_{ij}f_{jl} = \sum_{1 \leq j' \leq \vec{n}(x)} h_{ij'}h_{j'l}$
holds.
Therefore, the space of all such representations is the affine scheme
\begin{equation}   \label{A_V_ex2}
	\Rep(A , \vec{n}) \, = \, 
	\Spec \left(\, k \Big[ \, f_{ij} \, , \, g_{i'j'} \, , \, h_{i''j''} \, \Big] \, \Big/ \, 
	\Big(\sum_{1 \leq j \leq \vec{n}(y)} g_{ij}f_{jl} - \sum_{1 \leq j' \leq \vec{n}(x)} h_{ij'}h_{j'l} \Big) \, \right)
\end{equation}

In other words, to any presentation $A = k\langle Q \rangle /I$ by a quiver $Q$ with relations $I$, the space of representations $\Rep(A , \vec{n})$
is an affine scheme whose ring of functions is naturally identified with the commutative algebra
$A_{\vec{n}}$ obtained by formally inserting the subscript $(-)_{ij}$ on every generator and relation of $A$.

The representation scheme \eqref{A_V_ex2} allows us to construct the moduli space of 
finite dimensional representations of $A$.
Consider the group $\GL(\vec{n}) := \prod_{x\in \mathscr{O}} \, \GL(\vec{n}(x))$,
acting on the representation scheme $\Rep(A , \vec{n})$ by simultaneously 
conjugating each matrix generator of $A_{\vec{n}}$.
For example, consider the case \eqref{A_V_ex2}, the group $\GL(\vec{n})$ is simply $\GL(\vec{n}(x)) \times \GL(\vec{n}(y))$, so that 
an element $(g_x,g_y) \in \GL(\vec{n}(x)) \times \GL(\vec{n}(y))$ acts on $\Rep(A , \vec{n})$ by 
\begin{equation}   \label{GL_V_action_eq1}
	\Bigl( \, \bigl(  f_{ij}  \bigr) , \bigl(  g_{ij}  \bigr) , 
	\bigl(  h_{ij}  \bigr)    \, \Bigr)
	\, \mapsto \,
	\Bigl( \, 
	g_y \bigl(  f_{ij}  \bigr) g_x^{-1} , 
	g_x \bigl(  g_{ij}  \bigr) g_y^{-1}, 
	g_x \bigl(  h_{ij}  \bigr) g_x^{-1}    
	\, \Bigr)
\end{equation}

It is clear that two representations of $A$ are isomorphic if and only if they are conjugate under this action.
Accordingly, we define the \emph{moduli stack of representations} of $A$ with dimension vector $\vec{n}$ to be the quotient stack
\begin{equation}  \label{intro_Rep_stack}
\mathpzc{Rep}(A ; \vec{n}) \, := \, [ \, \Rep(A , \vec{n}) \, / \, \GL(\vec{n}) \, ] 
\end{equation}

This procedure of sticking subscript $(-)_{ij}$ to each generator and relation of a $k$-category works equally well when applied to a small dg category.
The result is then a commutative dg algebra, which would play the role of the ring of functions on the representation scheme.
This allows one to consider the derived version of \eqref{intro_Rep_stack} by first resolving $A$ by a cofibrant dg category $\cA$, and applying the same procedure.

For example, consider the dg category
\begin{equation}  \label{cA_V_ex1}
	\cA \, = \, k
	\left \langle
	\begin{tikzcd}
		x \ar[rr, bend left = 20, "f"]  \ar[loop left, distance=3em, "h"]
		\ar[loop above, "t"]
		& &  y \ar[ll, bend left = 20, "g"]
	\end{tikzcd} 
	\right \rangle
	, \qquad \qquad 
	d(t) =  gf - h^2
\end{equation}
where the generating arrow $t$ has cohomological degree $-1$, with differential  $d(t) = gf - h^2$. One can show that the canonical map $\cA \ronto A$ is a quasi-isomorphism, so that $\cA$ is a cofibrant resolution of $A$.

Thus, we take the commutative dg algebra
\begin{equation}  \label{cA_V_ex2}
	\cA_{\vec{n}} \, = \, 
	k \big[ \, f_{ij} \, , \, g_{i'j'} \, , \, h_{i''j''} \, , \, 
	t_{i''',j'''} \, \big] \, ,
	\qquad \quad
	d(t_{il}) \, = \, \sum_{1 \leq j \leq \vec{n}(y)} g_{ij}f_{jl} - \sum_{1 \leq j' \leq \vec{n}(x)} h_{ij'}h_{j'l}  
\end{equation}
and take the \emph{derived representation scheme} $\DRep(A; \vec{n}) := \Spec(\cA_{\vec{n}})$, which is now a derived scheme, and its quotient 
\begin{equation}   \label{intro_DRep_stack}
	 \mathpzc{DRep}(A ; \vec{n}) \, := \, [ \, \DRep(A , \vec{n}) \, / \, \GL(\vec{n}) \, ] 
\end{equation}
which is a derived stack, called the \emph{derived moduli stack of representations} of $\cA$. This construction is formalized in Section \ref{Rep_subsec}  and \ref{DRep_subsec}, where we perform the construction in a more coordinate-free way, and interpret the association $\cA \mapsto \cA_{\vec{n}}$ as a left adjoint functor, following \cite{BKR13}. We also establish in that section the relation between our derived moduli stack $\mathpzc{DRep}(A ; \vec{n})$ and the derived moduli stack $\cM_{A}$ of pseudo-perfect modules constructed in \cite{TV07}. Namely, we show in Theorem \ref{DRep_moduli_thm} that the open substack $\cM_{A}^{[0,0]} \subset \cM_{A}$ parametrizing modules of Tor-amplitude concentrated in $[0,0]$ is equivalent to the disjoint union
\begin{equation}   \label{intro_tor_ampl_equiv}
	\cM_{A}^{[0,0]} \, \simeq \,  
	\coprod_{\vec{n} : \scO \ra \bN} 
	\mathpzc{DRep}(A ; \vec{n})
\end{equation}

\subsection{Moduli space of representations}  \label{Rep_subsec}


We will study the moduli space of representations of a small $k$-category $A$. In the above, we have already seen how to proceed in a concrete example \eqref{A_V_ex1}, where we formally inserting the subscript $(-)_{ij}$ on every generator and relation of $A$, as in \eqref{A_V_ex2}.

This construction admits a conceptual description in terms of an adjoint functor.
For any collection $V = \{V_x\}_{x\in \mathscr{O}}$ of finite dimensional vector spaces, consider the endomorphism $k$-category $\Endcom(V)$ whose set of objects is $\mathscr{O}$, and whose morphism spaces are 
$\Endcom(V)(x,y) = \Hom_k(V_x, V_y)$, with the obvious composition maps.
Choosing a basis of each of the vector spaces $V_x$,
we may identify the collection $V$ as $V_x = k^{\vec{n}(x)}$
for some dimension vector $\vec{n}$. Then, a representation of $A$ with dimension vector $\vec{n}$ is the same as a $k$-linear functor $A \ra \Endcom(V)$ that fixes each object $x\in \mathscr{O}$.

Rename the commutative algebra $A_{\vec{n}}$ appearing in \eqref{A_V_ex2} as $A_V$, then the discussion that led to the consideration of \eqref{A_V_ex2} can be rephrased as saying that there exists a canonical bijection
\[
\Hom_{\Cat_{\scO,k}} (A , \Endcom(V)) \, \cong \, 
\Hom_{\CAlg_k}( A_V , k )
\]

This specifies the $k$-points of the scheme $\Rep(A,V) := \Spec(A_V)$ in terms of representations of $A$.
In fact, the above adjunction holds more generally. Namely, we have
\begin{equation}   \label{rep_adj2}
	\Hom_{\Cat_{\scO,k}} (A , \Endcom(V) \otimes B) \, \cong \, \Hom_{\CAlg_k}( A_V , B )
\end{equation}
for any commutative algebra $B \in \CAlg_k$. 
By Yoneda lemma, this adjunction completely determines the commutative algebra 
$A_V$ in a coordinate-free way, {\it i.e.}, without specifying a basis of $V$.
Our present discussion can then be summarized by the following

\bpp   \label{rep_adj_prop}
There exists an adjoint pair of functors 
\begin{equation*}
	\begin{tikzcd}
		(-)_V \, : \,  \Cat_{\scO,k} \ar[r, shift left = 0.3em]  
		&\CAlg_k \, : \,  \Endcom(V) \otimes -   \ar[l, shift left = 0.3em]
	\end{tikzcd}
\end{equation*}
\epp

We call the left adjoint functor $(-)_V$ the \emph{representation functor}, and the affine scheme $\Spec(A_V)$ the \emph{representation 
	scheme} of the $k$-category $A \in \Cat_k/\mathscr{O}$ with target $V = \{V_x\}_{x\in \scO}$.

Consider the unit map of this adjunction:
\begin{equation*} 
	\Phi \, : \, A \ra \Endcom(V) \otimes A_V
\end{equation*}
which is a representation of $A$ with coefficients in $A_V$,
known as the \emph{universal representation}.
More concretely, take any morphism $f \in A(x,y)$, together with $v \in V_x$ and $\lambda \in V_y^*$, then the triple $(\lambda,f,v)$ 
defines a function on the space $\Rep(A,V)$ of representations: any point $\rho \in \Rep(A,V)$ acts on $(\lambda,f,v)$ to give a value $\lambda(\rho(f)(v))$. Thus, the triple defines an element $\,_{\lambda} f_v \in A_V$. In other words, there is a map
\begin{equation}  \label{sticking_map_rep}
	V_y^* \otimes A(x,y) \otimes V_x \raq A_V
\end{equation}
which is in fact the map dual to the universal representation
\[
A(x,y) \raq V_y \otimes A_V \otimes V_x^* \, = \, \Endcom(V)(x,y) \otimes A_V
\]
%
%
If we choose basis $\{ v_j \} \subset V_x$ and $\{w_i\} \subset V_y$, with dual basis $\{ v_j^* \} \subset V_x^*$ and $\{w_i^*\} \subset V_y^*$,
then one can write the triple $(w_i^*, f,  v_j)$, considered as 
an element in $A_V$, as $f_{ij} \in A_V^{\ab}$.
%
Then the universal representation simply sends $f$ to the matrix $\bigl( f_{ij} \bigr) \in \bM_{\vec{n}(y) , \vec{n}(x)}(A_V) \cong \Endcom(V)(x,y) \otimes A_V$. 
Thus, one can view the universal representation as a formal way of saying that one can stick subscripts $(-)_{ij}$ to every morphism of $A$.

%
%

As in \eqref{GL_V_action_eq1}, the algebraic group $\GL_V := \prod_{x\in \mathscr{O}} \, \GL_{V_x}$
acts on the representation scheme $\Spec(A_V)$. 
More conceptually, for each $B \in \CAlg_k$, the group $\GL_V(B)$ acts on 
the $k$-category $\Endcom(V) \otimes B$ appearing on the left hand side of \eqref{rep_adj2}, which therefore induces an action on the right hand side of \eqref{rep_adj2}, functorial in $B$.
By Yoneda lemma, this gives an algebraic group action of $\GL_V$ on the affine scheme $\Spec(A_V)$.


The \emph{moduli stack of representations} of $A$ in $V$ is then 
defined as the quotient stack
\[
\mathpzc{Rep}(A ; V) \, := \, [ \, \Spec(A_V) \, / \, \GL_V \, ] 
\]
which is the stackification of the presheaf of groupoids
\[
[ \, \Spec(A_V) \, / \, \GL_V \, ]^{{\rm pre}} 
\, : \, \CAlg_k \raq {\rm Groupoid}
\]
represented by the algebraic action groupoid 
$\GL_V \ltimes \Spec(A_V)$.
In view of the adjunction \eqref{rep_adj2}, this presheaf sends 
$B \in \CAlg_k$ to the groupoid of $k$-linear functors
$T_B : A \ra \Mod(B)$ such that each $T_B(x) \in \Mod(B)$
is the \emph{free} module $V_x \otimes B$.
One can relax the last criterion by allowing $T_B(x)$ to be isomorphic, instead of equal, to $V_x \otimes B$.
This does not affect the equivalence type of the groupoid, so that 
we have
\begin{equation}  \label{Rep_prestack}
	[ \, \Spec(A_V) \, / \, \GL(V) \, ]^{{\rm pre}} \, (B)
	\, \simeq \,
	\Fun_k^{\circ} \bigl( \, A \,, \, \Mod^{{\rm free}}(B) \, \bigr)\Big|_{\, {\rm rank}(T(x)) \,=\, \dim(V_x)}
\end{equation}
where we have denoted $\Fun^{\circ}_k(C,D)$ the groupoid of 
$k$-linear functors from $C$ to $D$ and natural isomorphisms between functors, while the subscript to the end refers to the restriction 
to the full subcategory consisting of those functors 
$T: A \ra \Mod^{{\rm free}}(B)$ such that ${\rm rank}(T(x)) \,=\, \dim(V_x)$ for all $x\in \mathscr{O}$.

The stacification $[ \, \Spec(A_V) \, / \, \GL(V) \, ]$
of this presheaf of groupoid can then be described similarly
as the groupoid of functors 
to \emph{projective}, instead of free, $B$-modules, 
with prescribed ranks:
\begin{equation}  \label{Rep_stack}
	[ \, \Spec(A_V) \, / \, \GL(V) \, ] \, (B)
	\, \simeq \,
	\Fun_k^{\circ} \bigl( \, A \,, \, \Mod^{{\rm proj}}(B) \, \bigr)\Big|_{\, {\rm rank}(T(x)) \,=\, \dim(V_x)}
\end{equation}

\subsection{Derived representation schemes}  \label{DRep_subsec}

We extend the discussion in Section \ref{Rep_subsec} to dg categories.
Consider the category $\dgCatnOk$ whose objects are small (cohomologically) non-positively graded
dg categories with object set $\mathscr{O}$, and whose morphisms are dg functors that 
fix each object $x\in \mathscr{O}$.
(If $\mathscr{O}$ is finite, this is equivalent to the category of dg algebras 
over the semisimple ring $R = \bigoplus_{x\in \mathscr{O}} k$.)
Proposition \ref{rep_adj_prop} extends directly to this setting:

\bpp   \label{DRep_adj_prop}
There exists an adjoint pair of functors 
\begin{equation}  \label{DRep_adj1}
	\begin{tikzcd}
		(-)_V \, : \,  \dgCatnOk \ar[r, shift left = 0.3em]  
		&\CDGAn_k \, : \,  \Endcom(V) \otimes -   \ar[l, shift left = 0.3em]
	\end{tikzcd}
\end{equation}
\epp

As in the previous subsection, if we choose a basis for each $V_x$
to identify it with $k^{\vec{n}(x)}$, then the 
commutative dg algebra $\cA_V \in \CDGAn_k$ is simply obtained by 
formally inserting the subscripts $(-)_{ij}$ to morphisms of the dg category $\cA \in \dgCatnOk$. For example,
\eqref{cA_V_ex2} is the result of applying $(-)_V$ to 
\eqref{cA_V_ex1}.

The categories on both sides of \eqref{DRep_adj1} have canonical model structures where fibrations and weak equivalences are precisely 
those maps that induce fibrations and weak equivalences in all the 
relevant cochain complexes%
\footnote{
	When the object set $\mathscr{O}$ is finite, one can regard $\dgCatnOk$
	as the category of dg algebras over the semisimple algebra $R = k\mathscr{O}$,
	so that this falls precisely in the framework of \cite{BKR13}
	where one can find a detailed discussion.
	If the object set $\mathscr{O}$ is infinite, similar model structure exists,
	see, {\it e.g.}, \cite[Remark 4.6]{Tab10}.
}.
The functor $B \mapsto \Endcom(V) \otimes B$ clearly preserves 
weak equivalences and fibrations, which therefore make \eqref{DRep_adj1} a Quillen adjunction,
so that we have an induced adjunction at the level of homotopy categories
\begin{equation}  \label{DRep_adj2}
	\begin{tikzcd}
		{\bm L}(-)_V \, : \,  \Ho(\dgCatnOk) \ar[r, shift left = 0.3em]  
		&\Ho(\CDGAn_k) \, : \,  \Endcom(V) \otimes -   \ar[l, shift left = 0.3em]
	\end{tikzcd}
\end{equation}

For any $\cA \in \dgCatnOk$, one therefore has the 
(homotopy type of) commutative DG algebra ${\bm L}(\cA)_V$,
to which we can take $\Spec(-)$ (simply a formal symbol), which gives us an equivalence type 
of derived affine scheme, which we call the \emph{derived representation scheme}%
\footnote{We are using a slightly different notation than \cite{BKR13}, wherein the ``derived representation scheme'' 
	referred to the commutative dg algebra ${\bm L}(\cA)_V$ itself.}
of $\cA$ in $V$, denoted as
\begin{equation*}
	\DRep(\cA ; V) \, := \, \Spec( {\bm L}(\cA)_V )
\end{equation*}


As in the non-dg case, the unit map of the adjunction gives a universal representation
\begin{equation}  \label{univ_rep2}
	\Phi  \, : \,  \cA \ra \Endcom(V) \otimes \cA_V
\end{equation}
which formalizes the procedure of sticking subscripts $(-)_{ij}$.
By taking a cofibrant resolution $\cQ \rsa \cA$, one may also consider its derived version.

Like in the non-dg case, one can use the derived representation scheme to construct the derived moduli space of representations of $\cA$.
Indeed, for any cofibrant resolution $\cQ \rsa \cA$ in $\dgCatnOk$,
the commutative dg algebra $\cQ_V$ admits a $\GL_V$-action 
by simultaneous conjugation, completely analogous to \eqref{GL_V_action_eq1}.
This in turn allows one to take the derived stack
\begin{equation}  \label{stack_DRep_def}
	\mathpzc{DRep}(\cA ; V) \, := \, [ \, \Spec(\cQ_V) \, / \, \GL_V \, ]
\end{equation} 
which we call the \emph{derived moduli stack of representations} of 
$\cA$ into $V$.

\vspace{0.2cm}

We will now give a description of the moduli functor of the quotient stack $\mathpzc{DRep}(\cA ; V)$.
This will allow us to identify it with similar moduli stacks that have appeared in \cite{TV07}.
Recall that the notion of derived stacks can be formalized by considering simplicial presheaves on the category $\dAff = (\CDGAn_k)^{\op}$.
In the present case, the quotient stack $[ \, \Spec(\cQ_V) \, / \, \GL_V \, ]$ is the stackification of the simplicial presheaf
\begin{equation}   \label{DRep_prestack_1}
	[ \, \Spec(\cQ_V) \, / \, \GL_V \, ]^{{\rm pre}} \, : \, 
	\CDGAn_k \ra \Setdel, \qquad 
	B \mapsto 
	\big| \, \Map_{\dAff}( \, \Spec(B), N(\GL_V \ltimes \Spec(\cQ_V)) \, )    \, \big|
\end{equation}
where we have denoted by $N(\GL_V \ltimes \Spec(\cQ_V))$ 
the Segal groupoid of derived affine schemes corresponding to the 
algebraic action of $\GL_V$ on $\Spec(\cQ_V)$, so that 
the (levelwise) homotopy mapping space $\Map_{\dAff}( \, \Spec(B) \, , \, N( \GL_V \ltimes \Spec(\cQ_V) ) \, )$
is a simplicial space, of which we take the homotopy colimit $|-|$, which may be explicitly taken as the diagonal ({\it cf.} \cite[Corollary 18.7.7]{Hir03}).
The following result is a direct analogue to \eqref{Rep_prestack} and \eqref{Rep_stack}:

\bthm  \label{DRep_moduli_thm}
The simplicial presheaf $[ \, \Spec(\cQ_V) \, / \, \GL_V \, ]^{{\rm pre}}$ 
has an equivalent description
\begin{equation}  \label{quot_prestack_dgcat}
	[ \, \Spec(\cQ_V) \, / \, \GL_V \, ]^{{\rm pre}} (B)
	\, \simeq \,
	\Map_{\dgCat_k}  \bigl( \, \cQ \, , \, \Moddg^{{\rm free},0}(B)  \, \bigr) 
	\Big|_{\, {\rm rank}(T(x)) \,=\, \dim(V_x)}
\end{equation}
as the simplicial presheaf that assigns every $B \in \CDGAn_k$ 
the above homotopy mapping space taken in the model category 
$\dgCat_k$,
restricted to the connected components that consists of 
representations into free dg modules of specified ranks over $B$.

As a result, the stackification $[ \, \Spec(\cQ_V) \, / \, \GL_V \, ]$ of the simplicial presheaf $[ \, \Spec(\cQ_V) \, / \, \GL_V \, ]^{{\rm pre}}$ has an equivalent description
\[
[ \, \Spec(\cQ_V) \, / \, \GL_V \, ] (B)
\, \simeq \,
\Map_{\dgCat_k}  \bigl( \, \cQ \, , \, \Moddg^{{\rm proj},0}(B)  \, \bigr) 
\Big|_{\, {\rm rank}(T(x)) \,=\, \dim(V_x)}
\]
so that we have
\[
\cM_{\cA}^{[0,0]} \, \simeq \,  
\coprod_{\dim(V) :  \scO \ra \bN} 
\mathpzc{DRep}(\cA ; V)
\]
where $\cM_{\cA}^{[0,0]}$ is the derived moduli stack of pseudo-perfect dg modules on $\cA$ with Tor-amplitude contained in $[0,0]$ (see \cite{TV07} for definition, and \cite[Proposition 2.22.6]{TV07} for a characterization).
\ethm

To prove the first part, we will relate the homotopy pullback diagrams
\begin{equation} \label{BG_quot_pullback}
	\begin{tikzcd}
(\Spec \, \cQ_V)(B) \ar[r] \ar[d] &  {[} \, \Spec(\cQ_V) \, / \, \GL_V \, {]} ^{{\rm pre}}(B) \ar[d] \\
* \ar[r] & {[} \, * \, / \, \GL_V \, {]}^{{\rm pre}}(B) 
	\end{tikzcd}
\end{equation}
and
\begin{equation} \label{Map_dgcat_pullback}
	\begin{tikzcd}
		\Map_{k\scO \downarrow \dgCat_k}  \bigl( \, \cQ \, , \, \Moddg^{{\rm free},0}(B)  \, \bigr) \ar[r] \ar[d] &  \Map_{\dgCat_k}  \bigl( \, \cQ \, , \, \Moddg^{{\rm free},0}(B)  \, \bigr) \, =: \, X \ar[d] \\
		* \ar[r,"V_B"] & \Map_{\dgCat_k}  \bigl( \, k\scO \, , \, \Moddg^{{\rm free},0}(B)  \, \bigr)  \, =: \, Z
	\end{tikzcd}
\end{equation}
where $V_B$ is the point given by the functor $V_B : k\scO \ra \Moddg^{{\rm free},0}(B)$, that sends $x \in \scO$ to $V_x \otimes B$.

Consider the Quillen adjunction
\begin{equation}  \label{iota_R_Qadj}
	\begin{tikzcd}
		\iota \, : \, \dgCatnOk
		\ar[r, shift left]
		&  ( k\scO \downarrow \dgCatn_k ) \, : \, \mathfrak{R}
		\ar[l, shift left]
	\end{tikzcd}
\end{equation}
where $\iota$ is the obvious functor, while $\mathfrak{R}$ associates to $(F : k\scO \ra \cA)$ the dg category $\mathfrak{R}(F)$ on $\mathscr{O}$ such that $\mathfrak{R}(F)(x,y) := \cA(F(x),F(y))$.
Notice that there is an \emph{isomorphism} of dg categories
\begin{equation}  \label{R_VB_EndV}
	\mathfrak{R} \,( \, V_B : k\scO \ra \Moddg^{{\rm free},0}(B) \, ) \, \cong \, \Endcom(V) \otimes B
\end{equation}

The mapping spaces in \eqref{Map_dgcat_pullback} may, and will, be taken in $\dgCatn_k$ instead of $\dgCat_k$.
Then, the top left term has an alternative description:
\begin{equation}  \label{under_kO_V_adj}
	\begin{split}
	\Map_{k\scO \downarrow \dgCatn_k}  \bigl( \cQ , \Moddg^{{\rm free},0}(B)  \bigr) \, &\simeq \,
	\Map_{\dgCatnOk}  \bigl( \cQ , \Endcom(V) \otimes B) \\
	\, &\simeq \, \Map_{\CDGAn_k}(\cQ_V, B )
	\end{split}
\end{equation}
as a combination of the Quillen adjunctions \eqref{DRep_adj1} and \eqref{iota_R_Qadj}.

Denote by $Z_V$ the connected component of $Z$ that contains the point $V_B$, and by $X_V \subset X$ its preimage. 
Thus, $X_V$ is precisely the space appearing on the right hand side of \eqref{quot_prestack_dgcat}.
The homotopy equivalence \eqref{quot_prestack_dgcat} should then be intuitively clear. Namely, since $Z_V$ is connected, the space $X_V$ is equivalent to the total space of the Segal groupoid of the action of the based loop space of $Z$ at $V$ on the homotopy fiber. We have just seen that the homotopy fiber is equivalent to $\Map_{\CDGAn_k}(\cQ_V, B)$. Moreover, the based loop space of $Z$ at $V$ is the homotopy automorphism group of $V_B$ in $\Fun( k\scO , \Moddg^{{\rm free},0}(B)) \simeq  \Moddg^{{\rm free},0}(B)^{\times \scO}$, and hence is given by $\Map_{\CDGAn_k}(\cO(\GL_V),B)$ (see, {\it e.g.}, \cite{TV08}).
It should also be intuitively clear that the action of this based loop space on the fiber corresponds to the action of $\GL_V$ on $\Spec\, \cQ_V$.
However, instead of trying to make these identifications precise and functorial, we have found it easier to directly establish explicit simplicial sets that represent the mapping spaces in \eqref{BG_quot_pullback} and \eqref{Map_dgcat_pullback}, and show that they are homotopy equivalent (in fact, isomorphic). We now introduce the setup for that purpose.

\vspace{0.2cm}

Let $I$ be a small category. Then a \emph{Quillen system} on $I$ is an assignment $\cM_{\bullet}$ that assigns a model category $\cM_i$ for each $i \in \Ob(I)$, 
and a left Quillen functor $\varphi_* : \cM_i \ra \cM_j$ for each $\varphi \in \Hom_I(i,j)$, whose right adjoint will be denoted as $\varphi^* : \cM_j \ra \cM_i$. The system is required to be functorial in $\varphi$ up to a coherent system of natural isomorphisms, so that $\cM_{\bullet}$ is a pseudo-functor $I \ra {\rm CAT}$ to the $2$-category of (not necessarily small) categories.

Given a Quillen system $\cM_{\bullet}$ on $I$, then one can consider its category of sections $\Gamma(I,\cM_{\bullet})$. Namely, an object of $\Gamma(I,\cM_{\bullet})$ consists of a collection of object $X_i \in \cM_i$, one for each $i \in I$, as well as maps $\varphi_*(X_i) \ra X_j$ for each $\varphi \in \Hom_I(i,j)$, satisfying an obvious functoriality condition (involving the pseudo-functoriality of $\cM_{\bullet}$). 
Given $X_{\bullet}, Y_{\bullet} \in \Gamma(I,\cM_{\bullet})$, a morphism $F_{\bullet}:X_{\bullet}\ra Y_{\bullet}$ consists of a collection of maps $F_i \in \Hom_{\cM_i}(X_i,Y_i)$ commuting with the structure maps of $X_{\bullet}$ and $Y_{\bullet}$.
Notice that if $\cM_{\bullet}$ is the constant Quillen system on a model category $\cM$, then $\Gamma(I,\cM_{\bullet})$ is simply $\cM^I$. Our goal now is to generalize the result about the existence of the Reedy model structure on $\cM^I$ to the case of $\Gamma(I,\cM_{\bullet})$.

Thus, let $I$ be a Reedy category with the distinguished degree increasing (resp. decreasing) subcategory denoted by $I^+$ (resp. $I^-$), then for any $i \in \Ob(I)$ with $\deg(i) = \alpha$, define its latching and matching objects by taking the (co)limits in $\cM_i$
\begin{equation*}
	L_i X \, := \, \underset{(j,\varphi)\in I^+_{<\alpha}\downarrow i}{\colim} \, \varphi_*(X_j)
	\qquad \text{and} \qquad 
	M_i X \, := \, \underset{(j,\varphi)\in i \downarrow I^-_{<\alpha}}{\rm lim} \, \varphi^*(X_j)
\end{equation*}
As in the case of constant $\cM$, we have

\blm
For any given $X_{\bullet} \in \Gamma(I_{<\alpha},\cM_{\bullet})$, an extension of $X_{\bullet}$ to $\widetilde{X}_{\bullet} \in \Gamma(I_{\leq \alpha},\cM_{\bullet})$ is uniquely determined by a factorization
$L_i X \ra \widetilde{X}_i \ra M_iX$ for each $\deg(i)=\alpha$. 

Moreover, given $\widetilde{X}_{\bullet}, \widetilde{Y}_{\bullet} \in \Gamma(I_{\leq \alpha},\cM_{\bullet})$, and a given map $F_{\bullet} : X_{\bullet} \ra Y_{\bullet}$ on their restrictions to $\Gamma(I_{<\alpha},\cM_{\bullet})$, then an extension of $F_{\bullet}$ to $\widetilde{F}_{\bullet} : \widetilde{X}_{\bullet} \ra \widetilde{Y}_{\bullet}$ in $\Gamma(I_{\leq \alpha},\cM_{\bullet})$ is uniquely determined by maps $\widetilde{F}_i : \widetilde{X}_i \ra \widetilde{Y}_i$ in $\cM_i$, one for each $\deg(i)=\alpha$, that makes the following diagram commute:
\begin{equation*}
	\begin{tikzcd}
		L_iX \ar[d, "F"] \ar[r] & \widetilde{X}_i \ar[d, "\widetilde{F}_i"] \ar[r]  & M_iX \ar[d, "F"] \\
		L_iY \ar[r]  & \widetilde{Y}_i \ar[r]  & M_iY 
	\end{tikzcd}
\end{equation*}
\elm

\bpf
See, {\it e.g.}, \cite[Remark 5.2.3]{Hov99} for the case of constant $\cM$. The same proof carries through to non-constant Quillen systems. Alterntively, one can apply the result for constant $\cM$ to the Grothendieck construction associated to $\cM_{\bullet}$. 
\epf

Given a map $F_{\bullet}:X_{\bullet}\ra Y_{\bullet}$ in $\Gamma(I,\cM_{\bullet})$, the relative latching maps at $i \in I$ and the relative matching maps at $i \in I$ are respectively the following maps (both taken in $\cM_i$):
\begin{equation*}
	L_i Y \amalg_{L_i X} X_i \ra Y_i \qquad \text{and} \qquad X_i \ra Y_i \times_{M_i Y} M_i X
\end{equation*}

\blm
If each of the relative latching map is a (trivial) cofibration, then for each $i \in I$, the map $X_i \ra Y_i$ and the map $L_iX \ra L_iY$ are both (trivial) cofibration. The dual statements for relative matching map are also true.
\elm

\bpf
See {\it e.g.}, \cite[Corollary 5.1.5]{Hov99} for the case of constant $\cM$. The case for non-constant $\cM_{\bullet}$ also follows. Namely, consider the functor
\begin{equation*}
(I^+_{<\alpha}\downarrow i) \raq \cM_i \, , \qquad 	(j,\varphi) \, \mapsto \, \varphi_* X_j
\end{equation*}
and similarly for $Y$. Then $F$ determines a natural transformation between functors from the direct category $(I^+_{<\alpha}\downarrow i)$ to the model category $\cM_i$. One can check that the relative latching maps are still (trivial) cofibration, so we may directly apply \cite[Corollary 5.1.5]{Hov99} to conclude that $L_iX \ra L_iY$ is a (trivial) cofibration. The statement for $X_i \ra Y_i$ also follows because it can be written as the composition $X_i \ra L_i Y \amalg_{L_i X} X_i \ra Y_i$, both are which are (trivial) cofibrations.
\epf

As in the proofs of \cite[Theorem 5.1.3, 5.2.5]{Hov99}, the following result follows from the above two lemmae:
\bthm  \label{Reedy_thm}
There exists a model structure on $\Gamma(I,\cM_{\bullet})$ where a map $F_{\bullet}:X_{\bullet}\ra Y_{\bullet}$ is a 
\begin{enumerate}
	\item weak equivalence if each $F_i : X_i \ra Y_i$ is a weak equivalence in $\cM_i$.
	\item (trivial) cofibration if each relative latching map $L_i Y \amalg_{L_i X} X_i \ra Y_i$ is a (trivial) cofibration in $\cM_i$.
	\item (trivial) fibration if each relative matching map $X_i \ra Y_i \times_{M_i Y} M_i X$ is a (trivial) fibration in $\cM_i$.
\end{enumerate} 
\ethm

We will call this model structure the \emph{Reedy model structure}. The (trivial) cofibrations and (trivial) fibrations in this model category are called \emph{Reedy (trivial) cofibrations} and \emph{Reedy (trivial) fibrations}.

Now we apply this general setup to our case. 
For any map of sets $\varphi : \scO \ra \scO'$, notice that there is a Quillen adjunction
\begin{equation}  \label{dgcat_obj_map_adj}
\varphi_* \, : \,	\dgCatnOk \quad \substack{\longrightarrow \\ \longleftarrow} \quad  \dgCatn_{\scO',k} \, : \, \varphi^*
\end{equation}
where $\varphi^*$ is defined by $(\varphi^* \cA)(x,y) := \cA(\varphi(x),\varphi(y))$, and $\varphi_*(\cA)$ can be defined by sending a presentation of $\cA$ in terms of generators and relations to the same presentation but with source and target of the generating arrows and relations modified by $\varphi$.

%

Let $\scO$ be a given set. Consider the cosimplicial set $\widetilde{\scO}^{\bullet} : \Delta \ra \Set$ given by $\widetilde{\scO}^{n} := \coprod_{i \in [n]} \scO = \scO^{\amalg (n+1)}$. We will consider the Quillen system $\cM^{\bullet}$ on $\Delta$ given by $\cM^n := \dgCatn_{\widetilde{\scO}^n,k}$, with transition functors given by \eqref{dgcat_obj_map_adj}. 
Let $I^n$ be the category of $n+1$ objects with exactly one morphism between any two objects. We may regard $I^n$ as the category obtained by inverting every arrow in the poset category $[n]$. As such, we have $I^{\bullet} : \Delta \ra \Cat$. For any given $\cQ \in \dgCatnOk$, we then have an object
\begin{equation} \label{A_otimes_kI}
	\cQ \otimes k[I^{\bullet}] \in \Gamma(\Delta, \cM^{\bullet})
\end{equation}
where we recall that the tensor product of two small dg categories is defined by $\Ob(\cA \otimes \cB) = \Ob(\cA) \times \Ob(\cB)$ and $(\cA \otimes \cB)((x,y),(x',y')):= \cA(x,x') \otimes \cB(y,y')$. 

The tensor product \eqref{A_otimes_kI} has an alternative description as a pushout in $\dgCatnOk$:
\begin{equation}  \label{Q_otimes_kI_pushout}
		\cQ \otimes k[I^{n}] \, \cong \,  \cQ \amalg_{k\scO} (k\scO \otimes k[I^{n}])
\end{equation}
where the map $k\scO \ra k\scO \otimes k[I^{n}]$ is obtained by the inclusion $k \rinto k[I^{n}]$ into the $0$-th vertex.

Let $I$ and $J$ be Reedy categories, and $\cM_{\bullet}$ a Quillen system on $I$, then $\cM_{\bullet}^J$ is still a Quillen system on $I$, where each $\cM_i^J$ is endowed with the Reedy model structure. Thus, $\Gamma(I,\cM_{\bullet}^J)$ has a Reedy model structure by Theorem \ref{Reedy_thm}. Alternatively, we may rewrite this as $\Gamma(I,\cM_{\bullet}^J) \simeq \Gamma(I,\cM_{\bullet})^J$. The right hand side also has a model structure obtained by first applying Theorem \ref{Reedy_thm} to give a model structure on $\Gamma(I,\cM_{\bullet})$, which then induces a Reedy model structure on $(-)^J$. Clearly, these two model structures are the same, and will be called the \emph{double Reedy structure}. We will call $i \in I$ the horizontal degree and $j \in J$ the vertical degree.

We apply this to the case $I=J=\Delta$, and $\cM^n := \dgCatn_{\widetilde{\scO}^n,k}$ as above. From now on, fix a cofibrant $\cQ \in \dgCatnOk$. Think of \eqref{A_otimes_kI} as an object $(\cQ \otimes k[I^{\bullet}])^{vc} \in \Gamma(\Delta, (\cM^{\bullet})^{\Delta})$ that is constant in the vertical degree.
Resolve the morphism $(k\scO \otimes k[I^{\bullet}])^{vc} \ra (\cQ \otimes k[I^{\bullet}])^{vc}$ in $\Gamma(\Delta, (\cM^{\bullet})^{\Delta})$ under the (double Reedy) model structure of the last paragraph. {\it i.e.,} choose a commutative diagram
\begin{equation}  \label{scC_resolve}
	\begin{tikzcd}
		\scC_0 \ar[r] \ar[d] & \scC \ar[d] \\
		(k\scO \otimes k[I^{\bullet}])^{vc} \ar[r]
		& (\cQ \otimes k[I^{\bullet}])^{vc}
	\end{tikzcd}
\end{equation}
such that $\scC_0$ is cofibrant, $\scC_0 \ra \scC$ is a cofibration, and the vertical maps are weak equivalences. Notice that at the horizontal degree $0$, we have $k\scO \otimes k[I^{0}] = k\scO$, which is the initial object in $\cM^0$, so that it doesn't have to be resolved. We require that $\scC_0$ does not change this. {\it i.e.,} we require
\begin{equation}  \label{scC_0_triv}
	\scC_0^{0,j} \,= \, k\scO \qquad \text{for all } [j]\in \Delta
\end{equation}

\blm  \label{scC_cofib}
For each fixed vertical degree $j$, consider $\scC_0^{\bullet,j}$ and $\scC^{\bullet,j}$, and regard them as objects in $(\dgCatn_k)^{\Delta}$, endowed with the Reedy model structure. Then $\scC_0^{\bullet,j}$ is cofibrant and $\scC_0^{\bullet,j} \ra \scC^{\bullet,j}$ is a cofibration.
\elm

\bpf
The obvious (forgetful) functor
\begin{equation*}
\iota \, : \, \Gamma(\Delta, \dgCatn_{\widetilde{\scO}^{\bullet},k}) \raq (\dgCatn_k)^{\Delta}
\end{equation*}
sends latching objects to latching objects. {\it i.e.,} for any $\cA^{\bullet} \in \Gamma(\Delta, \dgCatn_{\widetilde{\scO}^{\bullet},k})$ and any $n \geq 1$, we have $\iota_n(L_n( \cA^{\bullet} )) = L_n(\iota(\cA^{\bullet}))$, where $\iota_n : \dgCatn_{\widetilde{\scO}^{n},k} \ra \dgCatn_k$ is the forgetful functor.
Therefore, $\iota$ sends Reedy cofibrant objects (resp. Reedy cofibrations) into Reedy cofibrant objects (resp. Reedy cofibrations). 
Indeed, one verifies that, after applying $\iota$, the latching map (resp. relative latching map) at $[n] \in \Delta$ is a cofibration. For $n \geq 1$, it follows from the above-mentioned fact that $\iota_n(L_n( \cA^{\bullet} )) = L_n(\iota(\cA^{\bullet}))$. The trivial case $n=0$ can be independently checked.
\epf

Let $X \ra Z$ be the map of bisimplicial sets
\begin{equation*}
	X_{i,j} := \Hom_{\dgCatn_k}( \scC^{i,j}, \Moddg^{{\rm free},0}(B) )
	\raq \Hom_{\dgCatn_k}( \scC_0^{i,j}, \Moddg^{{\rm free},0}(B) ) =: Z_{i,j}
\end{equation*}

By Lemma \ref{scC_cofib}, for each fixed vertical degree $j$, the map $X_{\bullet,j} \ra Z_{\bullet,j}$ is then an explicit model for the right hand side of \eqref{Map_dgcat_pullback} (see, {\it e.g.,} \cite[Definition 17.1.1]{Hir03}).
Moreover, $Z_{\bullet,j}$ is a Kan complex and $X_{\bullet,j} \ra Z_{\bullet,j}$ is a Kan fibration (see, {\it e.g.,} \cite[Corollary 16.5.3, 16.5.4]{Hir03}).
The relation with \eqref{BG_quot_pullback} is obtained by forming a certain first Eilenberg subcomplex. More precisely, we recall the following result (see, {\it e.g.,} \cite[Theorem 8.4]{May92}):

\blm
Let $Z$ be a Kan complex and $z \in Z_0$ be a vertex. Let $Z_{(z)}$ be the first Eilenberg subcomplex of $Z$ at $z$, {\it i.e.,} $Z_{(z)}$ consists of simplices of $Z$ whose vertices are $z$. Then the inclusion $Z_{(z)} \rinto Z$ is weakly equivalent to the connected component of $Z$ that contains $z$. 
\elm

We apply this Lemma to $Z_{\bullet,j}$. By our earlier choice \eqref{scC_0_triv}, the set $Z_{0,j}$ consists of an $\scO$-tuple of free dg modules over $B$. Let $V \in Z_{0,j}$ be the vertex corresponding the $\scO$-tuple $x \mapsto V_x \otimes B$. Let $Z_{0,j}^{(V)}$ be the corresponding first Eilenberg subcomplex, and let $X_{0,j}^{(V)} \subset X_{0,j}$ be its preimage. By our choice \eqref{scC_0_triv}, a dg functor $F: \scC^{i,j} \ra \Moddg^{{\rm free},0}(B)$ is in $X_{i,j}^{(V)}$ if and only if its effects on objects are given by $\widetilde{V}^i \otimes B$, where $\widetilde{V}^i$ is the collection of vector spaces $\widetilde{V}^i : \widetilde{\scO}^i = \scO^{i+1} \ra \scO \xra{V} \Vect_k$. The same is true for $Z_{i,j}^{(V)}$. In other words, we have
\begin{equation*}
	X_{i,j}^{(V)} = \Hom_{k\widetilde{\scO}^i \downarrow \dgCatn_k}( \scC^{i,j}, \Moddg^{{\rm free},0}(B) )
	\raq \Hom_{k\widetilde{\scO}^i \downarrow \dgCatn_k}( \scC_0^{i,j}, \Moddg^{{\rm free},0}(B) ) = Z_{i,j}^{(V)}
\end{equation*}

We may then apply the adjunctions \eqref{under_kO_V_adj} (at the level of Hom sets instead of mapping spaces), and get the description
\begin{equation}  \label{XV_ZV_ij}
	X_{i,j}^{(V)} = \Hom_{\CDGAn_k}( (\scC^{i,j})_{\widetilde{V}^i}, B )
	\raq \Hom_{\CDGAn_k}( (\scC_0^{i,j})_{\widetilde{V}^i}, B ) = Z_{i,j}^{(V)}
\end{equation}

We claim that this map of bisimplicial sets serves as an explicit model for the second column of \eqref{BG_quot_pullback}. To see this, we apply $(-)_{\widetilde{V}^{\bullet}}$ to \eqref{scC_resolve}. 
On the bottom row we get precisely the map $[\Spec \, \cQ_V / \GL_V] \ra [*/\GL_V]$ as a simplicial derived affine scheme. {\it i.e.,} we have
\begin{equation*}
	(k\scO \otimes k[I^{i}])_{\widetilde{V}^{i}} \, \cong \, \cO(\GL_V)^{\otimes i} \raq \cQ_V \otimes \cO(\GL_V)^{\otimes i} \, \cong \, (\cQ \otimes k[I^{i}])_{\widetilde{V}^{i}}
\end{equation*}
with the same cosimplicial maps as the one from the action groupoids. Indeed, this follows from the description \eqref{Q_otimes_kI_pushout}.

Notice that, for each fixed $i$, the objects $\scC^{i,\bullet}$ and $(\scC_0)^{i,\bullet}$ in $(\dgCatn_{\widetilde{\scO}^i,k})^{\Delta}$ are Reedy cofibrant, and hence so are the objects $(\scC^{i,\bullet})_{\widetilde{V}^i}$ and $(\scC^{i,\bullet}_0)_{\widetilde{V}^i}$ in $(\CDGAn_k)^{\Delta}$, since $(-)_{\widetilde{V}^i}$ is left Quillen.
Thus, to see that \eqref{XV_ZV_ij} represent the second column of \eqref{BG_quot_pullback}, it suffices to show that, for each $i,j$, the maps $(\scC^{i,j})_{\widetilde{V}^i} \ra (\cQ \otimes k[I^{i}])_{\widetilde{V}^{i}}$ and $(\scC^{i,j}_0)_{\widetilde{V}^i} \ra (k\scO \otimes k[I^{i}])_{\widetilde{V}^{i}}$ are quasi-isomorphisms of cdga's.
By construction, $\scC^{i,j} \ra \cQ \otimes k[I^{i}]$ and $\scC^{i,j}_0 \ra k\scO \otimes k[I^{i}]$ are cofibrant resolutions in $\dgCatn_{\widetilde{\scO}^i,k}$. Given a left Quillen functor $F : \cM \ra \cM'$, we say that $X \in \cM$ is adapted under $F$ if the canonical map ${\bm L}F(X) \ra F(X)$ in ${\rm Ho}(\cM')$ is an isomorphism. Thus, we are left to show the following

\bpp
The objects $\cQ \otimes k[I^{i}]$ and $k\scO \otimes k[I^{i}]$ are adapted under the functor $(-)_{\widetilde{V}^i} : \dgCatn_{\widetilde{\scO}^i,k} \ra \CDGAn_k$.
\epp 

\bpf
We may assume that $\cQ$ is semi-free. Then $\cQ \otimes k[I^{i}]$ is a semi-free extension of $k\scO \otimes k[I^{i}]$. Since $\CDGAn_k$ is left proper, it suffices to show the adaptedness of $k\scO \otimes k[I^{i}]$. Moreover, since $k\scO \otimes k[I^{i}]$ is just the $\scO$-fold disjoint union of $k[I^{i}]$, it suffices to prove it for the case $\scO = *$. 

In \cite{BFR14}, it is proved that if $A$ is a formally smooth associative algebra, then it is adapted under the functor $(-)_V : \DGA^{\leq 0}_k \ra \CDGAn_k$ for any finite dimensional vector space $V$. The same proof, based on \cite[Proposition 1.3]{Kap01}, can also be applied to the formally smooth small $k$-categories $k[I^{i}]$.
%
\epf

\bpf[Proof of Theorem \ref{DRep_moduli_thm}]
To summarize our above discussion, we see that the bisimplicial sets $X_{i,j}^{(V)}$ and $Z_{i,j}^{(V)}$ have the following properties:
\begin{enumerate}
	\item For each fixed $j$, the map of simplicial sets $X_{\bullet,j}^{(V)} \ra Z_{\bullet,j}^{(V)}$ is an explicit model for the map
	\begin{equation*}
		\Map_{\dgCatn_k}  \bigl( \cQ , \Moddg^{{\rm free},0}(B)  \bigr)\Big|_{ {\rm rank}(x) = \dim(V_x)} \ra \Map_{\dgCatn_k}  \bigl( k\scO , \Moddg^{{\rm free},0}(B)  \bigr)\Big|_{{\rm rank}(x) = \dim(V_x)} 
	\end{equation*}
\item For each fixed $i$, the map of simplicial sets $X_{i,\bullet}^{(V)} \ra Z_{i,\bullet}^{(V)}$ is an explicit model for the map
\begin{equation*}
	\Map_{\CDGAn_k}( \cQ_V \otimes \cO(\GL_V)^{\otimes i} ,B) \raq  	\Map_{\CDGAn_k}( \cO(\GL_V)^{\otimes i} ,B)
\end{equation*}
respecting the structure maps for varying $i$ coming from the action groupoids $[\Spec \, \cQ_V / \GL_V] \ra [*/\GL_V]$.
\end{enumerate}
Thus, if we take the diagonal simplicial set of $X_{\bullet,\bullet}^{(V)}$, it is simultaneously a model for both sides of \eqref{quot_prestack_dgcat}. 
The same is true for $Z_{\bullet,\bullet}^{(V)}$, so that we have a commutative diagram
\begin{equation*}
	\begin{tikzcd}
			{[} \, \Spec(\cQ_V) \, / \, \GL_V \, {]} ^{{\rm pre}}(B) \ar[r, "\simeq"] \ar[d]
		&   \Map_{\dgCat_k}  \bigl( \cQ , \Moddg^{{\rm free},0}(B)  \bigr)\Big|_{ {\rm rank}(x) = \dim(V_x)} \ar[d]
		\\
		{[} \, * \, / \, \GL_V \, {]}^{{\rm pre}}(B)   \ar[r, "\simeq"]
		& \Map_{\dgCat_k}  \bigl( k\scO , \Moddg^{{\rm free},0}(B)  \bigr)\Big|_{{\rm rank}(x) = \dim(V_x)} 
	\end{tikzcd}
\end{equation*}
where the horizontal maps are homotopy equivalences (in fact, isomorphisms).

To show the second statement, notice that if we replace $\Moddg^{{\rm free},0}$ in the second column by $\Moddg^{{\rm proj},0}$, then both spaces are derived stacks (see \cite[Corollary 1.3.7.4]{TV08}). By the universal property of stackifications, there is an induced map from the stackification of the left column to this modified right column. The induced map on the bottom row is then $\prod_{x \in \scO} B\GL(V_x)(B) \ra \prod_{x \in \scO} \Vect(B)|_{{\rm rank} = \dim(V_x)}$, which is known to be an equivalence (see, {\it e.g.}, \cite[Lemma 2.2.6.1]{TV08}). The induced map on homotopy fiber is also a homotopy equivalence since one can repeat the argument \eqref{under_kO_V_adj} with $\Moddg^{{\rm free},0}$ replaced by $\Moddg^{{\rm proj},0}$. Hence, the map on total space is also a homotopy equivalence.
\epf

\section{Shifted symplectic and Poisson structures on moduli spaces of representations}  \label{sec_sympl_poiss}

In Section \ref{sec_FNAG}, we mentioned that there are two aspects of formal noncommutative algebraic geometry, which we called the ``ontological aspect'' and the ``phenomenological aspect''. In this section, we develop some tools that allow us to work with the ``phenomenological aspect''. We apply these techniques to prove our main results, Theorems \ref{CY_induce_thm} and \ref{PCY_induce_thm}.

\subsection{Linearization and trace}  \label{subsec_trace}

We continue to assume that $k$ is a field of characteristic $0$. 
Let $\cA$ be a non-positively graded small dg category with $\scO = \Ob(\cA)$, and $V = \{V_x\}_{x \in \scO}$ a collection of finite dimensional vector spaces. In Section \ref{DRep_subsec}, we have constructed $\cA_V \in \CDGAn_k$, with an action of $\GL_V$.

In Section \ref{sec_FNAG}, we postulated that bimodules is a noncommutative analogue of quasi-coherent sheaves, and we also postulated the principle \eqref{NC_prin_naturalize}. We will establish the followings, which justify both of these postulate by the Kontsevich-Rosenberg principle:
\begin{enumerate}
	\item There is a functor $(-)_V^{\ab} : \Mod(\cAe) \raq \Mod_{\GL_V}(\cA_V)$.
	\item For any $M \in \Mod(\cAe)$, there is a map of complexes $\Psi : M_{\natural} \ra (M_V^{\ab})^{GL_V}$.
\end{enumerate}
The functor $(-)_V^{\ab}$ was called the \emph{Van den Bergh functor} in \cite{BKR13}, since it first appeared in \cite{VdB08b}. We will simply call it the \emph{linearization functor} (the name is justified by \eqref{vdB_sym_alg} below).
The map $\Psi$ will be called the \emph{trace map}.

Recall that the universal representation is the dg functor \eqref{univ_rep2}.
It allows us to consider the dg category $\Endcom(V) \otimes \cA_V$ as a bimodule over $\cA$. Moreover, since $\cA_V$ acts as central coefficients in this dg category, the $\cA$-bimodule structure and the $\cA_V$-module structure on $\Endcom(V) \otimes \cA_V$ are compatible.
Therefore, for any $M \in \Mod(\cAe)$, the tensor product $M \otimes_{\cAe} (\Endcom(V) \otimes \cA_V)$
is a dg module over $\cA_V$.
The linearization functor is given by (we discuss the $\GL_V$-equivariant structure below)
\begin{equation}   \label{vdB_functor}
	(-)_V^{\ab} \, : \, \Mod(\cAe) \raq \Mod(\cA_V),
	\qquad \qquad 
	M \, \mapsto \, M \otimes_{\cAe} (\Endcom(V) \otimes \cA_V)
\end{equation}

As observed in \cite{BKR13}, this functor is a left-adjoint functor.
Indeed, for any dg module $L \in \Mod(\cA_V)$, one can consider the bimodule $\Endcom(V) \otimes L$ over the dg category $\Endcom(V) \otimes \cA_V$. Again by the universal representation, one can consider $\Endcom(V) \otimes L$ as a bimodule over $\cA$.
This gives a functor
$\Endcom(V) \otimes - : \Mod(\cA_V) \raq \Mod(\cAe)$. 
A sequence of standard adjunction then shows the following
\blm 
There exists an adjoint pair of functors 
\begin{equation}  \label{vdB_adj1}
	\begin{tikzcd}
		(-)_V^{\ab} \, : \,  \Mod(\cAe) \ar[r, shift left = 0.3em]  
		& \Mod(\cA_V) \, : \,  \Endcom(V) \otimes -   \ar[l, shift left = 0.3em]
	\end{tikzcd}
\end{equation}
\elm

It is useful to consider the unit of this adjunction,
which is a map 
\begin{equation}  \label{vdB_unit_map}
	\Phi_{\cA} \, : \, M \raq \Endcom(V) \otimes M_V^{\ab}
\end{equation}
of bimodules over $\cA$.
Analogous to the unit map of the representation functor, 
this unit map allows us to think of the linearization functor as 
a formal procedure of sticking subscripts $(-)_{ij}$ to elements 
in $M$.
More precisely, consider an element $\xi \in M(x,y)$.
Then for any vector $v \in V_x$ and covector $\lambda \in V_y^*$, 
the triple $(\xi , v, \lambda)$ can be regarded as an element
\[
\xi \otimes ( \, ( v \otimes \lambda ) \otimes 1 \,)
\in M \otimes_{\cAe} (\Endcom(V) \otimes \cA_V) \, = \, M_V^{\ab}
\] 
This defines a map 
\begin{equation}  \label{vdB_sticking_map}
	M(x,y) \otimes V_x \otimes V_y^* \raq M_V^{\ab}
\end{equation}
which dually corresponds to a map 
\[
\Phi_{\cA} \, : \, M(x,y) \raq M_V^{\ab} \otimes V_y \otimes V_x^*
\, \cong \, \Endcom(V)(x,y) \otimes M_V^{\ab}
\]
which is precisely the $(x,y)$-component of the unit map \eqref{vdB_unit_map}.

Thus, if we choose basis $\{ v_j \} \subset V_x$ and $\{w_i\} \subset V_y$, with dual basis $\{ v_j^* \} \subset V_x^*$ and $\{w_i^*\} \subset V_y^*$,
then one can write the triple $(\xi , v_j, w_i^*)$, considered as 
an element in $M_V^{\ab}$, as $\xi_{ij} \in M_V^{\ab}$.
In this way, elements in $M_V^{\ab}$ are simply obtained by sticking subscripts $(-)_{ij}$ to elements in $M$.
Moreover, the unit map \eqref{vdB_unit_map} simply sends $\xi$
to the matrix $\bigl( \xi_{ij} \bigr) \in \Endcom(V)(x,y)\otimes M_V^{\ab}$.

The following lemma allows one to compute $M^{\ab}_V$ from a presentation of $M$:
\blm  \label{linearization_pres_lemma}
\begin{enumerate}
	\item The functor \eqref{vdB_functor} preserves colimits.
	\item Fix a choice of basis $\{ v_i^x \} \subset V_x$ for each $x \in \scO$. Suppose $M$ is semi-free over a set $\{ \xi_{\alpha} \in M(x_{\alpha},y_{\alpha})\}_{\alpha \in I}$, then $M^{\ab}_V$ is semi-free over the set $\{ (\xi_{\alpha})_{ij} \}_{\alpha \in I}$ (see the paragraph that precedes this lemma for the notation).
\end{enumerate}
\elm

\bpf
(1) is obvious as \eqref{vdB_functor} is a left adjoint (see \eqref{vdB_adj1}). (2) follows directly from the definition $M^{\ab}_V = M \otimes_{\cAe} (\Endcom(V) \otimes \cA_V)$.
\epf


An immediate consequence of the adjunction \eqref{vdB_adj1} is the natural isomorphism
\begin{equation}  \label{vdB_sym_alg}
	\bigl( \, T_{\cA}(M) \, \bigr)_V \, \cong \, \Sym_{\cA_V}(M_V^{\ab})
\end{equation}
One can show this isomorphism by explicitly checking generators and relations of both sides using Lemma \ref{linearization_pres_lemma}, or by a series of adjunctions that show that both sides (co)represent the same functor on the category of commutative dg algebras under $\cA_V$.
In fact, by endowing the bimodule $M$ with an extra weight grading $1$, 
the above isomorphism completely characterizes the linearization functor: the dg module $M_V^{\ab}$ is simply the weight graded $1$ component of $\bigl( \, T_{\cA}(M) \, \bigr)_V$.

We remark that this isomorphism respects the unit map for the representation functor and for the linearization functor.
Namely, take the universal representation \eqref{univ_rep2} for 
the dg category $T_{\cA}(M)$. Combining with the isomorphism \eqref{vdB_sym_alg} gives a dg functor
\[
T_{\cA}(M) \raq \Endcom(V) \otimes \bigl( \, T_{\cA}(M) \, \bigr)_V
\, \stackrel{\eqref{vdB_sym_alg}}{\cong} \, \Endcom(V) \otimes \Sym_{\cA_V}(M_V^{\ab})
\]
This dg functor respects the weight gradings on both sides. 
Restricting to weight $1$ components, one obtains a map 
$M \ra \Endcom(V) \otimes M_V^{\ab}$. 
It is easy to see that it coincides with the unit map 
\eqref{vdB_unit_map} for the linearization functor.
Thus, the isomorphism \eqref{vdB_sym_alg} characterizes both the linearization functor as well as the unit map of it. 
More importantly, it means that, for $\xi \in M(x,y)$, the meaning of the notation 
$\xi_{ij} \in M_V^{\ab}$ is unambiguous.

The isomorphism \eqref{vdB_sym_alg} allows us to give a $\GL_V$-equivariance structure on $M_V^{\ab}$. Namely, recall that giving a $\GL_V$-equivariance structure on a dg module $L$ over $\cA_V$ is equivalent to giving a $\GL_V$ action on $\Sym_{\cA_V}(L)$ extending the given action on $\cA_V$. Thus, the natural $\GL_V$-action on $\bigl(  T_{\cA}(M)  \bigr)_V $ provides such an equivariance structure on $M_V^{\ab}$.
Thus, we may regard the linearization functor as a functor
\begin{equation}   \label{vdB_functor_2}
	(-)_V^{\ab} \, : \, \Mod(\cAe) \raq \Mod_{\GL_V}(\cA_V)
\end{equation}

We now construct the above-mentioned trace map, {\it i.e.,} a map of cochain complexes
\begin{equation}   \label{vdB_Psi_map}
	\Psi  :  M_{\natural} \raq (M_V^{\ab})^{\GL_V}
\end{equation}

To construct this map, take any $\xi \in M(x,x)$ that represents an element in $M_{\natural}$. 
Apply the unit map \eqref{vdB_unit_map} to this element $\xi$ 
to obtain an element $\Phi_{\cA}(\xi) \in \Hom_k(V_x,V_x) \otimes M_V^{\ab}$. Taking trace, one has $\Tr(\Phi(\xi)) \in M_V^{\ab}$. 
Now if $\eta \in M(x,y)$ and $f \in \cA(y,x)$, then we have
$\Tr(\Phi(\eta \cdot f)) = (-1)^{|f||\eta|}\Tr(\Phi(f \cdot \eta))$. Therefore 
the map descends to a map on $M_{\natural}$. 
It is clear that the image of this map is contained in the 
dg submodule $(M_V^{\ab})^{\GL_V} \subset M_V^{\ab}$,
which then defines the map \eqref{vdB_Psi_map}.

We now discuss some properties of the linearization functor and the trace map. We will focus on two aspects: (1) effects under the monoidal products; (2) effects under (multi)duals.

Notice that both the domain and target of \eqref{vdB_functor_2} have canonical monoidal structures. 

\bpp  \label{op_lax_VdB_prop}
There is a canonical op-lax monoidal structure on the linearization functor \eqref{vdB_Psi_map}. {\it i.e.,} there is a canonical map in $\Mod_{\GL_V}(\cA_V)$:
\begin{equation}  \label{op_lax_unit_VdB}
	u \,: \, \cA_V^{\ab} \raq \cA_V
\end{equation}
For any $M,N \in \Mod(\cAe)$, there is a canonical map in $\Mod_{\GL_V}(\cA_V)$:
\begin{equation}  \label{op_lax_VdB}
	\psi_{M,N} \, : \, (\, M \otimes_{\cA} N \, )_V^{\ab} \raq M_V^{\ab} \otimes_{\cA_V} N_V^{\ab}
\end{equation}
which together satisfies associativity and unitality in the usual way.
\epp

\bpf
By the adjunction \eqref{vdB_adj1}, it suffices to give a lax monoidal structure on $\Endcom(V) : \Mod(\cA_V) \ra \Mod(\cAe)$. For $L_1, L_2 \in \Mod(\cA_V)$, it is given by taking the composition in $\Endcom(V)$:
\begin{equation}  \label{EndV_lax}
	(\Endcom(V) \otimes L_1 ) \otimes_{\cA} (\Endcom(V) \otimes L_2 ) \raq \Endcom(V) \otimes (L_1 \otimes_{\cA_V} L_2)
\end{equation}
The unitality map of this lax monoidal structure is given by the universal representation \eqref{univ_rep2}, regarded as a map of $\cA$-bimodules.

Unravelling the definition, this means that \eqref{op_lax_unit_VdB} is the map that corresponds under the adjunction \eqref{vdB_adj1} to the universal representation \eqref{univ_rep2}, regarded as a map of $\cA$-bimodules.
Likewise, for $M,N \in \Mod(\cAe)$, the structure map \eqref{op_lax_VdB} is the map that corresponds under the adjunction \eqref{vdB_adj1} to the map
\begin{equation*}
	M \otimes_{\cA} N \xraq{\eqref{vdB_unit_map}} 
	(\Endcom(V) \otimes M_V^{\ab}) \otimes_{\cA} (\Endcom(V) \otimes N_V^{\ab})
	\xraq{\eqref{EndV_lax}} \Endcom(V) \otimes (M_V^{\ab} \otimes_{\cA_V} N_V^{\ab})
\end{equation*}
\epf

\bcor  \label{Psi_multitrace_lemma}
For any $M_1,\ldots,M_n \in \Mod(\cAe)$, there is a canonical map of cochain complexes
\begin{equation}  \label{Psi_multitrace_map}
	\Psi \, : \, (\, M_1 \otimes_{\cA} \ldots \otimes_{\cA} M_n \,)_{\natural} \raq 
	\bigl( \, (M_1)^{\ab}_V \otimes_{\cA_V} \ldots \otimes_{\cA_V} (M_n)^{\ab}_V \, \bigr)^{\GL_V}
\end{equation}

Fix a choice of basis $\{ v_i^x \} \subset V_x$ for each $x \in \scO$. Then \eqref{Psi_multitrace_map} is given by
\begin{equation}  \label{Psi_multitrace_formula}
	\Psi(\xi_1 \otimes \ldots \otimes \xi_n) 
	\, = \, \sum_{i_0,\ldots,i_{n-1}}  (\xi_1)_{i_0 i_1} \otimes \ldots \otimes (\xi_n)_{i_{n-1} i_n}
\end{equation}
where we write $i_n = i_0$.
\ecor

\bpf
The map \eqref{Psi_multitrace_map} is obtained by combining \eqref{vdB_Psi_map} with the op-lax monoidal structure \eqref{op_lax_VdB}. Unravelling the definitions, it is given by first applying \eqref{vdB_unit_map}, then composing in $\Endcom(V)$, and then taking trace. This translates to the formula \eqref{Psi_multitrace_formula}. 
\epf

The linearization functor and the trace map also respects multi-duals. Recall that the procedure of forming $M_1^{\vee} \otimes_{\cA} \ldots \otimes_{\cA} M_n^{\vee}$ from $M_1,\ldots,M_n \in \Mod(\cAe)$ is not well-behaved because there are pathologies arising from taking tensor products of duals. Instead, the ``correct'' version is $\MD(M_1,\ldots,M_n)$. The two coincide if $M_1,\ldots,M_n$ are projective of finite rank as graded bimodules (see Lemma \ref{MD_map_isom_1}). Likewise, the ``correct'' version of $(M_1^{\vee} \otimes_{\cA} \ldots \otimes_{\cA} M_n^{\vee})_{\natural}$ is $\MD_{\natural}(M_1,\ldots,M_n)$.

Similar pathologies also arise in the commutative case. Thus, for a commutative dg algebra $B$ and for $L_1,\ldots,L_n \in \Mod(B)$, instead of taking $L_1^{\vee} \otimes_B \ldots \otimes_B L_n^{\vee}$, one should take $(L_1\otimes_B \ldots \otimes_B L_n)^{\vee}$. The two coincide if $L_1,\ldots,L_n$ are projective of finite rank as graded modules.
With these considerations in mind, we have

\bpp  \label{linearization_multidual_prop}
For any $M_1,\ldots,M_n \in \Mod(\cAe)$, there is a canonical map in $\Mod(\cA_V)$:
\begin{equation}   \label{linearization_multidual_module}
	\psi^{\dagger} \, : \,  \MD(M_1,\ldots, M_n)_V^{\ab} \raq 
	\bigl( \, (M_1)_V^{\ab}  \otimes_{\cA_V} \ldots  \otimes_{\cA_V} (M_n)_V^{\ab} \, \bigr)^{\vee}
\end{equation}
and a canonical map of cochain complexes%
\footnote{If $G$ is an affine algebraic group acting on a commutative dg algebra $B$, and if $L \in \Mod_G(B)$ is a $G$-equivariant dg module, then in general the dg module $L^{\vee} = \Homcom_B(L,B)$ may not have a $G$-equivariant structure. However, there is a cosimplicial object $[n] \mapsto \Homcom_B(L,B \otimes \cO(G)^{\otimes n})$ in $\Mod(B)$, which encapsulate most features of a $G$-equivariant structure (see, {\it e.g.}, \cite[Appendix A]{Yeu4}). In particular, one can define $(L^{\vee})^G$ as the equalizer of the two coface maps $d^0, d^1 : [1] \ra [2]$. In particular, the target of \eqref{linearization_multidual_naturalized} is defined in this way.}:
\begin{equation}   \label{linearization_multidual_naturalized}
	\Psi^{\dagger} \, : \,  \MD_{\natural}(M_1,\ldots, M_n) \raq 
	\Bigl( \, \bigl( \, (M_1)_V^{\ab}  \otimes_{\cA_V} \ldots  \otimes_{\cA_V} (M_n)_V^{\ab} \, \bigr)^{\vee} \, \Bigr)^{\GL_V}
\end{equation}
%
\epp

\bpf
To define \eqref{linearization_multidual_module}, given a pre-map 
$F \in \Homcom_{(\cA^{\otimes n})^e}( \, M_1 \otimes \ldots \otimes M_n \, , \, {}_{y}(\cA^{\otimes n+1})_{x} \, )$, consider the composition
\begin{equation}   \label{dagger_F}
	\begin{split}
		(M_1)_V^{\ab}  \otimes \ldots  \otimes (M_n)_V^{\ab} \, &= \, 
		(M_1 \otimes_{\cAe} (\Endcom(V) \otimes \cA_V)) \otimes \ldots \otimes (M_n \otimes_{\cAe} (\Endcom(V) \otimes \cA_V)) \\
		\, & = \,  (M_1 \otimes \ldots \otimes M_n) \otimes_{(\cA^{\otimes n})^e} (\Endcom(V) \otimes \cA_V)^{\otimes n} \\
		\, &\xra{F} \,
		{}_{y}(\cA^{\otimes n+1})_{x} \otimes_{(\cA^{\otimes n})^e} (\Endcom(V) \otimes \cA_V)^{\otimes n} \\
		\, & \xra{\eqref{univ_rep2}} 
		\, {}_{y}((\Endcom(V) \otimes \cA_V)^{\otimes n+1})_{x} \otimes_{(\cA^{\otimes n})^e} (\Endcom(V) \otimes \cA_V)^{\otimes n} \\
		\, & \xra{(*)} \, {}_{y}(\Endcom(V) \otimes \cA_V)_{x}
	\end{split}
\end{equation}
where the map $(*)$ is obtained by taking the composition in the dg category $\Endcom(V) \otimes \cA_V$ (notice that the tensor product over $(\cA^{\otimes n})^e$ in the second-to-last line guarantees that the object placements form a composable sequence). 

The pre-map \eqref{dagger_F} descends to an $\cA_V$-linear map from $(M_1)_V^{\ab}  \otimes_{\cA_V} \ldots \otimes_{\cA_V} (M_n)_V^{\ab}$, and is therefore an element in $\bigl( \, (M_1)_V^{\ab}  \otimes_{\cA_V} \ldots  \otimes_{\cA_V} (M_n)_V^{\ab} \, \bigr)^{\vee} \otimes {}_{y}\Endcom(V)_{x}$. This gives a map of bimodules $\MD(M_1,\ldots,M_n) \ra \bigl( \, (M_1)_V^{\ab}  \otimes_{\cA_V} \ldots  \otimes_{\cA_V} (M_n)_V^{\ab} \, \bigr)^{\vee} \otimes \Endcom(V)$. The map \eqref{linearization_multidual_module} is defined to be the one it corresponds to under the adjunction \eqref{vdB_adj1}.

Similarly, to define \eqref{linearization_multidual_naturalized}, given a pre-map $F \in \Homcom_{(\cA^{\otimes n})^e}( \, M_1 \otimes \ldots \otimes M_n \, , \, {}_{\tau}(\cA^{\otimes n})_{\id} \, )$, consider the composition
\begin{equation}  \label{Psi_dagger_F_1}
	\begin{split}
		(M_1)_V^{\ab}  \otimes \ldots  \otimes (M_n)_V^{\ab} 
 \, & = \,  (M_1 \otimes \ldots \otimes M_n) \otimes_{(\cA^{\otimes n})^e} (\Endcom(V) \otimes \cA_V)^{\otimes n} \\
 \, &\xra{F} \,
 {}_{\tau}(\cA^{\otimes n})_{\id} \otimes_{(\cA^{\otimes n})^e} (\Endcom(V) \otimes \cA_V)^{\otimes n} \\
 \, & \xra{\eqref{univ_rep2}} 
 \, {}_{\tau}((\Endcom(V) \otimes \cA_V)^{\otimes n})_{\id} \otimes_{(\cA^{\otimes n})^e} (\Endcom(V) \otimes \cA_V)^{\otimes n} \\
 \, & \xra{(**)} \, \cA_V
	\end{split}
\end{equation}
where the map $(**)$ is obtained by taking the composition in the dg category $\Endcom(V) \otimes \cA_V$ and then taking trace (notice that the tensor product over  $(\cA^{\otimes n})^e$ in the second-to-last line guarantees that the object placements form a composable sequence, and whose composition lies in $\End(V_x) \otimes \cA_V$ for some $x \in \scO$, to which we can take the trace $\End(V_x) \ra k$).

The map \eqref{Psi_dagger_F_1} descends to a $\GL_V$-equivariant $\cA_V$-linear pre-map
\begin{equation*}
	(M_1)_V^{\ab}  \otimes_{\cA_V} \ldots  \otimes_{\cA_V} (M_n)_V^{\ab} 
	\raq \cA_V
\end{equation*}
which we define to be $\Psi^{\dagger}(F)$.
\epf

\bcor \label{VdB_dual_cor}
For any $M \in \Mod(\cAe)$, there is a canonical map in $\Mod(\cA_V)$
\begin{equation}  \label{psi_dagger_single}
	\psi^{\dagger} \, : \, (M^{\vee})^{\ab}_V \raq (M^{\ab}_V)^{\vee}
\end{equation}
which is an isomorphism if $M$ is projective of finite rank as a graded bimodule.
\ecor

\bpf
The map \eqref{psi_dagger_single} is defined to be \eqref{linearization_multidual_module} for $n = 1$. 
If $\lambda \in {}_{y}(M^{\vee})_{x}$ and $\xi \in {}_{y'}M_{x'}$, and if we write by Sweedler's notation
\begin{equation*}
	\langle \lambda , \xi \rangle \, = \, \langle \lambda , \xi \rangle' \otimes \langle \lambda , \xi \rangle'' \, \in \, {}_{y}\cA_{x'} \otimes {}_{y'}\cA_{x}
\end{equation*}
then we have 
\begin{equation}  \label{psi_dagger_single_formula}
\langle \, \psi^{\dagger}(\lambda_{ij}) \, , \, \xi_{pq} \, \rangle \, = \, 
\langle \lambda , \xi \rangle'_{iq} \cdot \langle \lambda , \xi \rangle''_{pj}
\end{equation}

If $M$ is semi-free over a finite set $\{  \xi_{\alpha} \}$, then by Lemma \ref{linearization_pres_lemma}, $(M^{\vee})^{\ab}_V$ is semi-free over $\{  ((\xi_{\alpha})^{\vee})_{ji} \}$, while $(M^{\ab}_V)^{\vee}$ is semi-free over $\{  ((\xi_{\alpha})_{ij})^{\vee} \}$. The formula \eqref{psi_dagger_single_formula} shows that \eqref{psi_dagger_single} sends $((\xi_{\alpha})^{\vee})_{ji}$ to $((\xi_{\alpha})_{ij})^{\vee}$, and hence is an isomorphism in this case.
The case when $M$ is projective of finite rank as a graded bimodule then follows by taking a retract (the statement is independent of the differential).
\epf

The canonical maps in Proposition \ref{op_lax_VdB_prop} and \ref{linearization_multidual_prop} satisfy some compatibility conditions. We will only make use of the following:

\blm  \label{Psi_trace_comm_lemma}
For any $M,N \in \Mod(\cAe)$, the following diagram is commutative:
\begin{equation}  \label{Psi_trace_comm_diag}
	\begin{tikzcd}
		(M \otimes_{\cA} N)_{\natural} \ar[rr] \ar[d, "\eqref{Psi_multitrace_map}"']
		 & & \Homcom_{\cAe}(M^{\vee},N) \ar[d, "(-)^{\ab}_V"] \\ 
		M^{\ab}_V \otimes_{\cA_V} N^{\ab}_V \ar[r]
		&  \Homcom_{\cA_V}((M^{\ab}_V)^{\vee}, N^{\ab}_V ) \ar[r, "\eqref{psi_dagger_single}^*"]
		& \Homcom_{\cA_V}((M^{\vee})^{\ab}_V, N^{\ab}_V ) 
	\end{tikzcd}
\end{equation}
where the vertical map $(-)^{\ab}_V$ on the right is the functoriality map on Hom complexes of the canonical dg enrichment of \eqref{vdB_functor}.
\elm

\bpf
Given $\omega \in (M \otimes_{\cA} N)_{\natural}$. For simplicity, assume that it is given by $\omega = \xi \otimes \eta$ (in general it is a finite sum of such).
By \eqref{Psi_multitrace_formula}, we have $\Psi(\omega) = \sum_{p,q} \xi_{pq} \otimes \eta_{qp}$. Hence, by \eqref{psi_dagger_single_formula}, the lower route of \eqref{Psi_trace_comm_diag} sends $\omega$ to the pre-map
\begin{equation*}
	\lambda_{ij} \, \mapsto \, \sum_{p,q} \langle \lambda , \xi \rangle'_{iq} \cdot \langle \lambda , \xi \rangle''_{pj} \cdot \eta_{qp} 
	\, = \, (\langle \lambda , \xi \rangle' \cdot \eta \cdot \langle \lambda , \xi \rangle'')_{ij} \, = \, (\omega^{\#}(\lambda))_{ij}
\end{equation*}
This is precisely the image of $\omega$ under the upper route of \eqref{Psi_trace_comm_diag}.
\epf

We now discuss the derived version of the linearization functor. Clearly, \eqref{vdB_adj1} is a Quillen adjunction, so that we have a derived functor
${\bm L}(-)^{\ab}_V : \cD(\cAe) \ra \cD(\cA_V)$. 
Denote by $(-)_{h\natural} : \cD(\cAe) \ra \cD(k)$ the derived version of $(-)_{\natural}$, and by $(-)^!$ the derived version of $(-)^{\vee}$ on $\cD(\cAe)$ and on $\cD(\cA_V)$, then we have the following analogue of Corollary \ref{VdB_dual_cor} and Lemma \ref{Psi_trace_comm_lemma}:

\bcor
\begin{enumerate}
	\item For any $M \in \cD(\cAe)$, there is a canonical map in $\cD(\cA_V)$
	\begin{equation}  \label{psi_dagger_single_derived}
		\psi^{\dagger} \, : \, {\bm L}(M^{!})^{\ab}_V \raq ({\bm L}(M)^{\ab}_V)^{!}
	\end{equation}
	which is an isomorphism if $M$ is perfect.
	\item For any $M,N \in \cD(\cAe)$, the following diagram is commutative:
	\begin{equation*}  
		\begin{tikzcd}
			(M \otimes^{{\bm L}}_{\cA} N)_{h\natural} \ar[rr] \ar[d, "\eqref{Psi_multitrace_map}"']
			& & \RHomcom_{\cAe}(M^{!},N) \ar[d, "{\bm L}(-)^{\ab}_V"] \\ 
			{\bm L}(M)^{\ab}_V \otimes^{{\bm L}}_{\cA_V} {\bm L}(N)^{\ab}_V \ar[r]
			&  \RHomcom_{\cA_V}(({\bm L}(M)^{\ab}_V)^{!}, {\bm L}(N)^{\ab}_V ) \ar[r, "\eqref{psi_dagger_single_derived}^*"]
			& \RHomcom_{\cA_V}({\bm L}(M^{!})^{\ab}_V, {\bm L}(N)^{\ab}_V ) 
		\end{tikzcd}
	\end{equation*}
\end{enumerate}
\ecor

\bcor  \label{VdB_nondeg_cor}
Suppose that $M,N \in \Mod(\cAe)$ are cofibrant, and $M$ is perfect. If an element $\omega \in Z^m((M \otimes_{\cA} N)_{\natural})$ is (left) non-degenerate in the sense that the induced map $\omega^{\#} : M^{\vee} \ra N$ in $\Mod(\cAe)$ is a quasi-isomorphism, then its image $\Psi(\omega) \in M^{\ab}_V \otimes_{\cA_V} N^{\ab}_V$ under \eqref{Psi_multitrace_map} is also (left) non-degenerate in the sense that the induced map $\Psi(\omega)^{\#} : (M^{\ab}_V)^{\vee} \ra N^{\ab}_V$ is a quasi-isomorphism.
\ecor

\subsection{Shifted symplectic and Poisson structures on global quotients}

Let $G$ be a reductive group acting on a non-positively graded commutative dg algebra $B$, we recall the main results of \cite{Yeu4}, which gives an explicit characterization of shifted symplectic and Poisson structures on the derived stack $X = [ Y / G ]$ for $Y = \Spec \, B$.
The reader may find more details in \cite{Yeu4}. 

In order to distinguish with the bimodule of noncommutative Kahler differential, we will denote by $\Omega^1_{\com}(B)$ the usual dg module of Kahler differentials of a commutative dg algebra $B$.
Denote by $\alpha : \Omega^1_{\com}(B) \ra \frg^* \otimes B$ the infinitesimal action, which is $G$-equivariant. Define $\Omega^1_{\Car}(Y/G) \in \Mod_G(B)$ by
\begin{equation}  \label{Omega_Car_def}
	\Omega^1_{\Car}(Y/G)[1] \, = \, \cone[ \, \Omega^1_{\com}(B) \xra{\alpha} \frg^* \otimes B  \, ]
\end{equation}
If $B$ is cofibrant, then this equivariant module represents the cotangent complex of $X = [Y/G]$.

Taking global sections in $\QCoh(X)$ corresponds to taking the $G$-invariants of a $G$-equivariant dg module. This gives a description of the global $n$-forms $\DR_{\Car}^n(Y/G)$ of $X$ as a cochain complex
\begin{equation}  \label{DRCar_1}
	\DR_{\Car}^n(Y/G)[-n] \, = \, (\Sym_B^n( \Omega^1_{\Car}(Y/G)[-1] ))^G \, = \, \bigoplus_{p+q = n} ( \, \Omega^p_{\com}(B)[-p] \otimes \Sym^q(\frg^*[-2]) \,)^G
\end{equation}

The differential on $\DR_{\Car}^n(Y/G)$ will be denoted by $\partial$.
The differential $D : \Omega^{\bullet}_{\com}(B) \ra \Omega^{\bullet + 1}_{\com}(B)$ induces a map $D' = D \otimes \id : \DR_{\Car}^{\bullet}(Y/G) \ra \DR_{\Car}^{\bullet+1}(Y/G)$. The following is a classical result (see, {\it e.g.,} \cite{Yeu4} for a proof):

\bpp  \label{Cartan_DR_prop}
The maps $D' : \DR_{\Car}^{n}(Y/G) \ra \DR_{\Car}^{n+1}(Y/G)$ is a map of cochain complex $\partial D' = D' \partial$, and satisfies $D' \circ D' = 0$. Hence $(\DR_{\Car}^{\bullet}(Y/G),\partial,D')$ is a bicomplex, called the \emph{Cartan-de Rham bicomplex}.
\epp

We take the direct product total complex
\begin{equation}
	\DR_{\Car}^{\cl}(Y/G) \,:= \, \Pi_{p \geq 0} \, \DR_{\Car}^p(Y/G)[-p] \, , \qquad \qquad \dtot = \partial + D'
\end{equation}
which comes with the \emph{Hodge filtration}
\begin{equation}
	F^r \, \DR_{\Car}^{\cl}(Y/G) \,:= \, \Pi_{p \geq r} \, \DR_{\Car}^p(Y/G)[-p] \, , \qquad \qquad \dtot = \partial + D'
\end{equation}

It is clear that any cocycle $\omega \in Z^m(\DR^2(Y/G))$ induces a map
\begin{equation}  \label{omega_sharp}
	\omega^{\sharp} \, : \, \Omega_{\Car}^1(Y/G)^{\vee}[-m] \raq \Omega_{\Car}^1(Y/G)
\end{equation}
of dg modules over $B$.

\bdf
We say that $B \in \CDGAn_k$ is \emph{almost cofibrant} if $\Omega^1_{\com}(B) \in \Mod(B)$ is cofibrant, and if the map $\bL_{B/k} \ra \Omega^1_{\com}(B)$ is a quasi-isomorphism.
\edf

\bdf  \label{shifted_sympl_def}
Suppose that $G$ is reductive and $B$ is almost cofibrant, then an \emph{$m$-shifted pre-symplectic structure in Cartan model} on $[Y/G]$ is a cocycle $\widetilde{\omega} \in Z^{m+2}F^2 \DR_{\Car}^{\cl}(Y/G)$.

Denote by $\omega \in Z^m (\DR^2(Y/G))$ the part of $\widetilde{\omega}$ in Hodge degree $2$ ({\it i.e.}, it is the image under the quotient by $F^3$). Then $\widetilde{\omega}$ is said to be \emph{non-degenerate} if the corresponding map \eqref{omega_sharp} is a quasi-isomorphism. 
In this case, we say that $\widetilde{\omega}$ is an \emph{$m$-shifted symplectic structure in Cartan model} on $[Y/G]$.
\edf

Now we discuss the dual picture of (shifted) polyvector fields and (shifted) Poisson structures. For any $B \in \CDGAn_k$, define the complex of \emph{$m$-shifted $p$-polyvector fields} by
\begin{equation*}
	\Pol^p(B,m) \, := \, \Homcom_B( \, \Sym^p_B(\Omega^1_{\com}(B)[m+1]) \, , \, B \,)
\end{equation*}

The Schouten-Nijenhuis bracket is defined by considering $\Pol^*(B,m)$ as a subcomplex
\begin{equation*}
	\Homcom_B( \, \Sym^p_B(\Omega^1_{\com}(B)[m+1]) \, , \, B \,) \, \subset \, \Homcom_k(  (B[m+1])^{\otimes p}  , B )
\end{equation*}
consisting of those maps that are a derivation in each variable, and is symmetric under the $S_n$-action.
The right hand side (shifted by $m+1$) has a dg pre-Lie algebra structure given by
\begin{equation}  \label{pre_Lie_shuffle_1}
	f * g \, := \, \sum_{i=1}^p \, \sum_{\sigma \in S^i_{p,q}} \,(f \circ_i g)^{\sigma}
\end{equation}
for $f \in \Homcom_k(  (B[m+1])^{\otimes p}  , B[m+1] )$ and $g \in \Homcom_k(  (B[m+1])^{\otimes q}  , B[m+1] )$, where 
\begin{equation*}
	S^i_{p,q} \,:= \, \Bigl\{ \, \sigma \in S_{p+q-1} \, | \, \substack{\sigma^{-1}(1)<\sigma^{-1}(2) <\ldots< \sigma^{-1}(i) < \sigma^{-1}(i+q) < \ldots < \sigma^{-1}(p+q-1) \\
		\text{and } \sigma^{-1}(i) < \sigma^{-1}(i+1) < \ldots < \sigma^{-1}(i+q-1)} \, \Bigr\}
\end{equation*}

Write $\langle n \rangle := \{1,\ldots,n\}$. For any subset $S \subset \langle p+q-1 \rangle$ of cardinality $|S|=q$, 
there exists a unique pair $(i,\sigma)$ where $1 \leq i \leq p$ and $\sigma \in S^i_{p,q}$, such that $\sigma(S) = \{i, i+1,\ldots,i+q-1\}$. (Here, $i$ is simply the smallest element of $S$). Conversely, the pair $(i,\sigma)$ uniquely determines $S = \sigma^{-1}(\{i, i+1,\ldots,i+q-1\})$.
For any $1 \leq j \leq p$, denote by $g_{S,j} : (B[m+1])^{\otimes p+q-1} \ra (B[m+1])^{\otimes p}$ the map that sends $(B[m+1])^{\otimes S}$ to the $j$-th copy of $B[m+1]$ via $g$ ($S$ inherits a total order from $\langle p+q-1 \rangle$, so that $(B[m+1])^{\otimes S}$ is canonically identified with $(B[m+1])^{\otimes q}$), and the identity map on the rest (in an order-preserving way on the tensor indexing). Write $g_S := g_{S,j}$ for $j = \min(S)$.
Then \eqref{pre_Lie_shuffle_1} can be rewritten as
\begin{equation} \label{pre_Lie_shuffle_2}
	f * g \, := \, \sum_{S \subset \langle p+q-1 \rangle, \, \, |S|=q} f \circ (g_{S})
\end{equation}

One can show that the associated Lie bracket $\{f,g\} := f * g - (-1)^{|f||g|} g * f$ preserves the subcomplex $\Pol^*(B,m)[m+1]$ (see, {\it e.g.}, \cite[Section 2]{Mel16}), which we take as the Schouten-Nijenhuis bracket $\{-,-\}$, so that $\Pol^*(B,m)[m+1]$ is a dg Lie algebra.

If $G$ is a reductive group acting on $B$, then we define
\begin{equation*}
	\begin{split}
		\Pol_{\Car}^n(Y/G,m) \, &:= \, \Homcom_B( \, \Sym^n_B( \Omega_{\Car}^1(Y/G)[m+1] ) \, , \, B \,)^G  \\
		\, &= \, \bigoplus_{p+q = n} ( \, \Pol^p(B,m) \otimes {\rm CoSym}^q(\frg[-m]) \,)^G
	\end{split}
\end{equation*}

Define the bracket
\begin{equation}  \label{SN_brac_Cartan}
	\{-,-\}' \, : \, \Pol_{\Car}^p(Y/G,m)[m+1] \otimes \Pol_{\Car}^q(Y/G,m)[m+1] \raq \Pol_{\Car}^{p+q-1}(Y/G,m)[m+1]
\end{equation}
by $\{-,-\}' = \{-,-\} \otimes \mu$, where $\{-,-\}$ is the Schouten-Nijenhuis bracket on $\Pol^{\bullet}(B,m)$ and $\mu$ is the shuffle product on ${\rm CoSym}^{\bullet}(\frg[-m])$.
Then we have (see \cite[Theorem 2.13]{Yeu4}):
\bthm
The bracket \eqref{SN_brac_Cartan} makes $\Pol_{\Car}^{\bullet}(Y/G,m)[m+1]$ into a dg Lie algebra.
\ethm

We write a formula for this bracket. Given $G \in \Pol_{\Car}^{q}(Y/G,m)[m+1]$, 
then for any subset $S \subset \langle p+q-1 \rangle$ of cardinality $|S|=q$ and any $1 \leq j \leq p$, denote by
\begin{equation*}
	\widetilde{G}_{S,j} \, : \, (\Omega_{\Car}^1(Y/G)[m+1])^{\otimes p+q-1} \raq (\Omega_{\Car}^1(Y/G)[m+1])^{\otimes p}
\end{equation*}
the map that sends $(\Omega_{\Car}^1(Y/G)[m+1])^{\otimes S} \xra{G} B[m+1] \xra{D} \Omega_{\Car}^1(Y/G)[m+1]$ to the $j$-th component, and the identity map on the rest (in an order-preserving way on the tensor indexing). Write $\widetilde{G}_{S} := \widetilde{G}_{S,j}$ for $j = \min(S)$. 


Then, for $F \in \Pol_{\Car}^{p}(Y/G,m)[m+1]$, $G \in \Pol_{\Car}^{q}(Y/G,m)[m+1]$, and $\omega \in (\Omega_{\Car}^1(Y/G)[m+1])^{\otimes p+q-1}$, if we define
\begin{equation}  \label{star_Pol_Car}
	(F * G)(\omega) \, := \, \sum_{S \subset \langle p+q-1 \rangle, \, \, |S|=q} \, F \bigl( \widetilde{G}_{S}(\omega) \bigr)  
\end{equation}
then we have%
\footnote{Notice that $(F * G)(\omega)$ is well-defined for $\omega \in (\Omega_{\Car}^1(Y/G)[m+1])^{\otimes p+q-1}$ where we take the tensor over $k$. When we take the commutator \eqref{bracket_Pol_Car}, it descends to a map from the tensor over $B$.} 
\begin{equation} \label{bracket_Pol_Car}
	\{F,G\}'(\omega) = (F * G)(\omega) - (-1)^{|F||G|} (G * F)(\omega)
\end{equation}

Indeed, the analogue of \eqref{bracket_Pol_Car} in the non-equivariant setting follows from \eqref{pre_Lie_shuffle_2}. As in \cite{Yeu4}, $\Pol_{\Car}^{*}(Y/G,m)[m+1]$ can be regarded as a dg Lie subalgebra of the dg Lie algebra $\Pol^{*}([Y/\frg],m)[m+1]$ of polyvectors on the Chevalley-Eilenberg cdga $[Y/\frg]$. The formula \eqref{bracket_Pol_Car} is then inherited from the analogous one for $[Y/\frg]$.

To define shifted Poisson structures, take the completion
\begin{equation*}
	\widehat{\Pol}_{\Car}(Y/G,m) \, := \, \prod_{p \geq 0} \, \Pol_{\Car}^p(Y/G,m)
\end{equation*}
which comes with the weight filtration
\begin{equation*}
	F^r\widehat{\Pol}_{\Car}(Y/G,m) \, := \, \prod_{p \geq r} \, \Pol_{\Car}^p(Y/G,m)
\end{equation*}

The Schouten-Nijenhuis-Cartan bracket makes $\widehat{\Pol}_{\Car}(Y/G,m)[m+1]$ a dg Lie algebra, and each of the subcomplexes $F^r\widehat{\Pol}_{\Car}(Y/G,m)[m+1]$ a dg Lie subalgebra.

\bdf  \label{shifted_Poiss_def}
Suppose that $G$ is reductive and $B$ is almost cofibrant, then an \emph{$m$-shifted Poisson structure in Cartan model} on $[Y/G]$ is a Maurer-Cartan element in the dg Lie algebra $F^2\widehat{\Pol}_{\Car}(Y/G,m)[m+1]$.
%
\edf

The main result of \cite{Yeu4} is the following

\bthm
The notions of shifted symplectic and Poisson structures in Definitions \ref{shifted_sympl_def} and \ref{shifted_Poiss_def} coincide with the ones in \cite{Pri17}.
\ethm

\subsection{Shifted symplectic and Poisson structures on moduli spaces of representations}

Let $\cA$ be a non-positively graded small dg category with object set $\scO$, and let $V = \{V_x\}_{x \in \scO}$ be a collection of finite dimensional vector spaces. Let $B = \cA_V$ with the given action by $G = \GL_V$, as in Section \ref{DRep_subsec}. We will consider shifted symplectic and Poisson structures on the global quotients $X = [Y/G]$ where $Y = \Spec \, B$. Our main results are Theorems \ref{CY_induce_thm} and \ref{PCY_induce_thm}.

\blm  \label{Omega_ab_V_lemma}
There is a canonical isomorphism $\Omega^1(\cA)^{\ab}_V \cong \Omega^1_{\com}(\cA_V)$.
\elm

\bpf
The map $D : \cA_V \ra \Omega^1_{\com}(\cA_V)$ induces a derivation 
\begin{equation*}
	\id \otimes D  \, : \, \Endcom(V) \otimes \cA_V \raq \Endcom(V) \otimes \Omega^1_{\com}(\cA_V)
\end{equation*}
of the bimodule $\Endcom(V) \otimes \Omega^1_{\com}(\cA_V)$ over the dg category $\Endcom(V) \otimes \cA_V$. 

Composing with the universal representation \eqref{univ_rep2}, we have a derivation $(\id \otimes D) \circ \Phi : \cA \ra \Endcom(V) \otimes \Omega^1_{\com}(\cA_V)$, which therefore corresponds to a map of bimodules $\Omega^1(\cA) \ra \Endcom(V) \otimes \Omega^1_{\com}(\cA_V)$. By the adjunction \eqref{vdB_adj1}, this corresponds to a map
\begin{equation}  \label{Omega_ab_V_map}
	\Omega^1(\cA)^{\ab}_V \raq \Omega^1_{\com}(\cA_V) \, , \qquad (Df)_{ij} \mapsto D(f_{ij})
\end{equation}

When $\cA$ is semi-free, say $\cA = T_{\scO}(Q)$, then $\Omega^1_{\com}(\cA_V)$ is semi-free over $\{ D(f_{ij}) \}_{f \in Q}$, while $\Omega^1(\cA)^{\ab}_V$ is semi-free over $\{ (Df)_{ij} \}_{f \in Q}$ by Lemma \ref{linearization_pres_lemma}(2). Thus \eqref{Omega_ab_V_map} is an isomorphism when $\cA$ is semi-free. In general, one can present $\cA$ as $\cA = T_{\scO}(Q)/I$. The corresponding presentations for $\Omega^1(\cA)^{\ab}_V$ (see Lemma \ref{linearization_pres_lemma}(1)) and for $\Omega^1_{\com}(\cA_V)$ get identified under \eqref{Omega_ab_V_map}, so that \eqref{Omega_ab_V_map} is an isomorphism in general.
\epf

\bcor  \label{cSA_ab_V_cor}
There is a canonical isomorphism $\cSA^{\ab}_V \cong \Omega^1_{\Car}(Y/G)[1]$.
\ecor

\bpf
Apply $(-)^{\ab}_V$ to the cone $\cSA = \cone[\, \Omega^1(\cA) \xra{\alpha} \cA \otimes_{\scO} \cA \,]$. By Lemma \ref{Omega_ab_V_lemma}, the first term becomes $\Omega^1_{\com}(\cA_V)$. By Lemma \ref{linearization_pres_lemma}(2), the second term $(\cA \otimes_{\scO} \cA)^{\ab}_V$ is given by $\bigoplus_{x \in \scO} \, \End(V_x) \otimes \cA_V = \mathfrak{gl}_V^* \otimes \cA_V$. Under these identifications, the map $\alpha^{\ab}_V$ is given by $\alpha^{\ab}_V(D(f_{ij})) = \sum_{p} f_{ip}E_{pj} - \sum_{q} E_{iq}f_{qj}$, where we regard $E_{ij}$ as elements in $\mathfrak{gl}_V^*$. This is precisely the infinitesimal action map that defines the cone \eqref{Omega_Car_def}. 
In other words, the desired isomorphism is given by
\begin{equation}  \label{cSA_ab_V_map}
	\cSA^{\ab}_V \xraq{\cong} \Omega^1_{\Car}(Y/G)[1] \, , \qquad (sDf)_{ij} \mapsto sD(f_{ij}) \, , \quad E_{ij} \mapsto E_{ij}
\end{equation}
\epf

Combining this with \eqref{Psi_multitrace_map}, we obtain a map
\begin{equation}  \label{Psi_cSA_1}
	\Psi \, : \, ( \, \cSA \otimes_{\cA} \stackrel{(r)}{\ldots} \otimes_{\cA} \cSA \,)_{\natural} \raq (\, \Omega^1_{\Car}(Y/G)[1] \otimes_{\cA_V} \stackrel{(r)}{\ldots} \otimes_{\cA_V} \Omega^1_{\Car}(Y/G)[1]  \,)^{\GL_V}
\end{equation}

The domain of \eqref{Psi_cSA_1} has a canonical $C_r$-action, while the target of \eqref{Psi_cSA_1} has a canonical $S_r$-action. The map \eqref{Psi_cSA_1} intertwines with these actions under the inclusion $C_r \subset S_r$. Passing to the coinvariants under these actions, we have an induced map
\begin{equation}  \label{Psi_cSA_2}
	\Psi \, : \, \scX^{(r)}(\cA) \raq \DR_{\Car}^r(Y/G)[r]
\end{equation}

Recall from \eqref{scX_mixed} that $\{\scX^{(r)}(\cA)\}_{r \geq 0}$ has the structure of an $\bN$-graded mixed structure. On the other hand, the fact that $\DR_{\Car}^{\bullet}(Y/G)$ is a bicomplex (see Proposition \ref{Cartan_DR_prop}) can be translated to the statement that $\{\DR_{\Car}^r(Y/G)[r]\}_{r \geq 0}$ has the structure of an $\bN$-graded mixed structure. We have
\bpp  \label{Psi_cSA_mixed_prop}
The maps \eqref{Psi_cSA_2} is a map of $\bN$-graded mixed complexes.
\epp

\bpf
An element of $\scX^{(r)}(\cA)$ is a finite sum of elements of the form
\begin{equation*}
	\Theta \, = \, f_0 \cdot \theta_1 \cdot f_1 \cdot \ldots \cdot f_{r-1} \cdot \theta_r \cdot f_r
\end{equation*}
where $f_i \in {}_{x_i}\cA_{y_{i+1}}$, and $\theta_i$ is either of the form $\theta_i = sDg_i$ for $g_i \in {}_{y_i}\cA_{x_{i}}$ or $\theta_i = E_{x_i}$ (in which case $x_i = y_i$), such that $y_{r+1} = x_0$. 

As we mentioned in Section \ref{sec_nc_DR}, the mixed map $B$ on $\{\scX^{(r)}(\cA)\}_{r \geq 0}$ is obtained by writing the letter $D$ in front of the expression, and simplify using the Leibniz rule. In particular, we see that both $\Psi(B(\Theta))$ and $sD'(\Psi(\Theta))$ are given by
\begin{equation*}
	 \sum_{p = 0}^r (-1)^{\spadesuit} \sum_{i_0,\ldots,i_r, j_1,\ldots,j_r} \,
	(f_0)_{i_0 j_1} \cdot (\theta_1)_{j_1 i_1} \cdot (f_1)_{i_1 j_2} \cdot \ldots \cdot (sDf_p)_{i_p j_{p+1}} \cdot \ldots \cdot (f_{r-1})_{i_{r-1} j_r} \cdot (\theta_r)_{j_r i_r} \cdot (f_r)_{i_r i_0} 
\end{equation*}
where the Koszul sign is given by $\spadesuit = |f_0| + |\theta_1| + \ldots + |f_{p-1}| + |\theta_p|$.
\epf

\bthm  \label{CY_induce_thm}
Let $\cA$ be smooth. Then any $n$-Calabi-Yau structure on $\cA$ induces a $(2-n)$-shifted symplectic structure on the derived moduli stack of representations $\mathpzc{DRep}(\cA ; V)$.
\ethm

\bpf
Assume without loss of generality that $\cA$ is cofibrant. By Proposition \ref{Psi_cSA_mixed_prop}, there is a map of complexes
\begin{equation*} 
	\Psi \, : \, F^2\scX^{\tot}(\cA) \raq F^2\DR_{\Car}^{\cl}(Y/G)
\end{equation*}
By Corollary \ref{VdB_nondeg_cor}, if $\widetilde{\omega} \in Z^{4-n}(F^2\scX^{\tot}(\cA))$ is non-degenerate, then so is $\Psi(\widetilde{\omega}) \in Z^{4-n}(F^2\DR_{\Car}^{\cl}(Y/G))$.
\epf

Now we consider the dual picture of polyvector fields. Combining Corollary \ref{cSA_ab_V_cor} with \eqref{linearization_multidual_naturalized}, we obtain a map
\begin{equation}  \label{Psi_dagger_MD_cSA}
	\Psi^{\dagger} \, : \,\MD^{(r)}_{\natural}( \cSA[m] ) \raq 
	\bigl( \,( \, \Omega^1_{\Car}(Y/G)[m+1] \otimes_{\cA_V} \stackrel{(r)}{\ldots} \otimes_{\cA_V} \Omega^1_{\Car}(Y/G)[m+1] \,)^{\vee}  \, \bigr)^{\GL_V}
\end{equation}

By taking a sum over cosets $S_r/C_r$, we obtain a map from the $C_r$-invariants to the $S_r$-invariants:
\begin{equation}  \label{Psi_Pol_map}
	\Psi^{\dagger}_{S_r/C_r} \, : \, \scP^{(r)}(\cA,m) \xraq{ \sum_{\sigma \in S_r/C_r} \sigma_* \circ \Psi^{\dagger} } \Pol_{\Car}^r(Y/G,m)
\end{equation}


\bthm  \label{polyvec_trace_dgla}
The maps \eqref{Psi_Pol_map} give a map of dg Lie algebras 
\begin{equation*}
\bigoplus_{r \geq 1} \, \scP^{(r)}(\cA,m)[m+1] \raq \bigoplus_{r \geq 1} \, \Pol_{\Car}^r(Y/G,m)[m+1]
\end{equation*}
\ethm

\bpf
Fix a choice of basis of each $V_x$. Consider the data $\vec{\theta} = (\theta_1,\ldots,\theta_p)$, where $\theta_v$ is either $sDf_v$ for some $f_v \in {}_{y_v}\cA_{x_v}$, or $E_{x_v}$ (in which case we write $y_v = x_v$). In the case $\theta_v = sDf_v$, we think of $f_v$ as part of the data for $\theta_v$.
Also, consider $\vec{i} = (i_1,\ldots,i_r)$ and $\vec{j} = (j_1,\ldots,j_r)$, where $i_v$ is an indexing element for the basis of $V_{y_v}$, and $j_v$ is an indexing element for the basis of $V_{x_v}$.

Associated to this data, define
\begin{equation*}
\omega(\vec{\theta},\vec{i},\vec{j},m) \, := \, s^{m}(\theta_1)_{i_1 j_1} \otimes \ldots \otimes s^{m}(\theta_p)_{i_p j_p} \, \in \, (\Omega^1_{\Car}(Y/G)[m+1])^{\otimes p}
\end{equation*}
where we write $(sDf_v)_{i_v j_v} := sD((f_v)_{i_v j_v})$. Here, the tensor product is taken over $k$. Denote by
\begin{equation*}
	\overline{\omega}(\vec{\theta},\vec{i},\vec{j},m) \, \in \, \Omega^1_{\Car}(Y/G)[m+1] \otimes_{\cA_V} \stackrel{(p)}{\ldots} \otimes_{\cA_V} \Omega^1_{\Car}(Y/G)[m+1]
\end{equation*}
its image in the tensor product over $\cA_V$.

Given $F \in \scP^{(p)}(\cA,m)$, we may apply it to $(\theta_1,\ldots,\theta_p)$ to obtain
\begin{equation*}
	F(\vec{\theta},m) \,:= \, F(s^{m}\theta_1 \otimes \ldots \otimes s^{m}\theta_p) \, \in \, {}_{y_p}\cA_{x_1} \otimes \ldots \otimes {}_{y_{p-1}}\cA_{x_p}
\end{equation*}
Alternatively, we may think of the data $(\theta_1,\ldots,\theta_p)$ as determining an $\scO$-colored disk together with input data, to which we apply $F$. In this way, $F(\vec{\theta},m)$ makes sense also for $F \in C_H^{(p)}(\cA,m)^{C_p}$.

Pairing with $\tau(\vec{i}) := (i_p,i_1,\ldots,i_{p-1})$ and $\vec{j}$, and then multiply within $\cA_V$, we obtain an element
\begin{equation*}
	F(\vec{\theta},m)_{(\tau(\vec{i}),\vec{j})} \, \in \, \cA_V
\end{equation*}

By definition, the map \eqref{Psi_dagger_MD_cSA} sends $F \in \scP^{(p)}(\cA,m)$ to the pre-map
\begin{equation*}
\Psi^{\dagger}(F) \, : \, \Omega^1_{\Car}(Y/G)[m+1] \otimes_{\cA_V} \stackrel{(p)}{\ldots} \otimes_{\cA_V} \Omega^1_{\Car}(Y/G)[m+1] \raq \cA_V \, , \qquad 
	\overline{\omega}(\vec{\theta},\vec{i},\vec{j},m) \mapsto F(\vec{\theta},m)_{(\tau(\vec{i}),\vec{j})}
\end{equation*}

Given data $(\vec{\theta},\vec{i},\vec{j})$ of length $p+q-1$, write $\omega := \omega(\vec{\theta},\vec{i},\vec{j},m)$. For any $\sigma \in S_{p+q-1}$, write  $\sigma(\omega) = \sigma(\omega)' \otimes \sigma(\omega)''$, where \begin{equation*}
	\begin{split}
		\sigma(\omega)' \, &:= \, s^{m}(\theta_{\sigma(1)})_{i_{\sigma(1)} j_{\sigma(1)}} \otimes \ldots \otimes s^{m}(\theta_{\sigma(p-1)})_{i_{\sigma(p-1)} j_{\sigma(p-1)}}  \\
		\sigma(\omega)'' \, &:= \, s^{m}(\theta_{\sigma(p)})_{i_{\sigma(p)} j_{\sigma(p)}} \otimes \ldots \otimes s^{m}(\theta_{\sigma(p+q-1)})_{i_{\sigma(p+q-1)} j_{\sigma(p+q-1)}}
	\end{split}
\end{equation*}

Given $F \in \scP^{(p)}(\cA,m)[m+1]$ and $G \in \scP^{(q)}(\cA,m)[m+1]$, then we have
\begin{equation}  \label{FoG_theta}
  \sum_{\sigma \in S_{p+q-1}/C_{p+q-1}} (F \circ G)(\sigma(\vec{\theta}), m)_{(\tau(\sigma(\vec{i})), \sigma(\vec{j}))}
  \, = \, \frac{1}{q} \sum_{\sigma \in S_{p+q-1}} (-1)^{\spadesuit} \,  \Psi^{\dagger}(F) \biggl( \sigma(\omega)' \otimes D \Bigl( \Psi^{\dagger}(G) \bigl( \sigma(\omega)'' \bigr) \Bigr) \biggr)
\end{equation}
with the Koszul sign given by $\spadesuit = |G||\sigma(\omega)'|$. The factor $\frac{1}{q}$ on the right hand side appears because the derivation property of $D(-)$ on the right hand side gives rise to $q$ terms, which are repetitive under the sum over $\sigma \in S_{p+q-1}$. For the left hand side, $F \circ G$ is a sum over the possible position of compositions (which are $C_{p+q-1}$ rotations of each other), so that the sum over $S_{p+q-1}/C_{p+q-1}$ becomes a sum over $S_{p+q-1}$ if we fix a choice of position of composition. Incorporating both of these, we see the equality \eqref{FoG_theta}.

On the other hand, we compute the product \eqref{star_Pol_Car} between $\Psi^{\dagger}_{S_p/C_p}(F)$ and $\Psi^{\dagger}_{S_q/C_q}(G)$:
\begin{equation}  \label{Psi_FG_star}
\bigl( \Psi^{\dagger}_{S_p/C_p}(F) * \Psi^{\dagger}_{S_q/C_q}(G) \bigr) (\omega) \, = \, 
\sum_{S \subset \langle p+q-1 \rangle, \, \, |S|=q} \, \frac{1}{pq}  \sum_{\sigma' \in S_p} \sum_{\sigma'' \in S_q} \Psi^{\dagger}(F)^{\sigma'} \Bigl( \, (\widetilde{\Psi^{\dagger}(G)^{\sigma''}})_{S}  ( \omega )  \,\Bigr) 
\end{equation}

Since $\Psi^{\dagger}(F)$ is $C_p$-invariant, the sum $\frac{1}{p}\sum_{\sigma' \in S_p}$ on the right hand side can be replaced by the sum $\sum_{\sigma' \in S_p , \,\, \sigma'(\min(S)) = p}$. 
Given $S'$ and $S''$ any totally ordered set of cardinality $r$, then any permutation $\sigma \in S_r$ determines a canonical bijection $\sigma : S' \xra{\cong} S''$, and vice versa. 
Given $\sigma' \in S_p$ such that $\sigma'(\min(S)) = p$, denote by $\sigma^{\prime \circ} :  \langle p  \rangle \setminus \{\min(S)\}  \xra{\cong} \langle p-1 \rangle $ its restriction, which therefore determines a bijection $\sigma^{\prime \circ} : \langle p+q-1 \rangle \setminus S \xra{\cong} \langle p-1 \rangle$. 
Together with $\sigma''$, this uniquely determines $\sigma \in S_{p+q-1}$ that sends $S$ to $\{ p,\ldots,p+q-1 \}$ via $\sigma''$, and is given by $\sigma^{\prime \circ}$ on $\langle p+q-1 \rangle \setminus S$. This gives a bijection 
\begin{equation*}
	\{ \,  (S,\sigma',\sigma'')  \, | \, S \subset \langle p+q-1 \rangle, \, |S|=q, \, \sigma' \in S_p, \, \sigma'' \in S_q, \, \sigma'(\min(S)) = p \, \} \, \cong \, S_{p+q-1}
\end{equation*}
Under this bijection, each term on the right hand side of \eqref{Psi_FG_star} becomes the corresponding term on the right hand side of \eqref{FoG_theta}. Therefore, \eqref{FoG_theta} and \eqref{Psi_FG_star} are equal.
\epf

As a corollary, we have

\bthm  \label{PCY_induce_thm}
Any $n$-pre-Calabi-Yau structure on $\cA$ induces a $(2-n)$-shifted Poisson structure on the derived moduli stack of representations $\mathpzc{DRep}(\cA ; V)$.
\ethm

\bpf
Assume without loss of generality that $\cA$ is cofibrant. Take $m = 2-n$. By Theorems \ref{ext_necklace_Lie_subalg} and \ref{polyvec_trace_dgla}, there are maps of weight graded dg Lie algebras
\begin{equation*}
	\bigoplus_{r \geq 1} \, C_H^{(r)}(\cA,m)^{C_r}[m+1] \xlintoq{\sim} \bigoplus_{r \geq 1} \, \scP^{(r)}(\cA,m)[m+1] \raq \bigoplus_{r \geq 1} \, \Pol_{\Car}^r(Y/G,m)[m+1]
\end{equation*}
where the left pointing map is a quasi-isomorphism. Take the associated maps on Maurer-Cartan spaces (defined either as the mapping space $\underline{{\rm MC}}(L) := \Map_{{\rm dgla}_k^{{\rm gr}}}(k(-2)[-1],L)$ in the model category of weight graded dg Lie algebras, or more explicitly as in \cite[Definition 1.5]{Pri17}).
\epf

\end{document}